# The Least Odd Partition Identity in Conjunction with Dirichlet Linear Progressions and the Reduced Residue System Modulo 18;
# A potential efficient deterministic P=NP solution for the P vs. NP problem of primality testing and factorization

Introducing the core lowest odd partition identity (LOPI) in conjunction with the reduced residue system Modulo (18), a deterministic system of Dirichlet binomial linear progressions for sequential prime/non-prime element generation that upon specific linear convolution partitions the proper non-prime subset matrices via generation of all non-trivial factors; a potentially new polynomial time efficient factorization/prime testing algorithm; ELMA-18.

Laurel L. McClure
laurellynn.mcclure@studio.unibo.it
Bologna, Italy

Abstract: In this work we introduce the numerical constant, LOPI, $N_{LOPI} \equiv LOPI\ (modulo\ 18) = Integer\ Sum,$ the lowest odd partition identity in conjunction with the reduced residue system Modulo 18, a complete disjoint covering residue system when considered in its whole set of residues of ,{0,1,2,3…m-1}. By convolution of two specific LOPI Dirichlet linear progressions per LOPI subset, the elements of the 6 reduced residue congruence class sets containing all possible prime and non-prime natural odd numbers congruent to Modulo18 are generated sequentially and uniquely. Through set generation and complete composite number subset generation based on the pattern of sequential factor generation for all composite numbers in each LOPI congruence class, multiplicity becomes a possible tool to define degree of primality.  We show that the lowest sum of the digits of an integer, even or odd, is a constant unrestricted partition identity contained each natural number and is equal to one of the residues in Mod 18, {1, 5, 7, 11, 13, 17}.  We term this constant the lowest odd/even partition identity or LOPI/LEPI.  When integers are viewed as partitions of the core LOPI/LEPI value, the infinite natural numbers represent the residues of Mod 18 as extended each time toward infinity by sequential cycles of 18.  LOPI mod 18 is a factor generation system that demonstrates why the primes >18 are not generated through convolution of specific Dirichlet linear progressions (LOPI$_1$ + 18(a))(LOPI$_2$ + 18(b)) into the LOPI non-prime subset matrices and reveals the pattern of the primes to be actually six cyclic patterns. By then querying the multiplicity in these 3-4 LOPI specific non-prime matrices we suggest another definition of prime and non-prime based on this multiplicity, m≥0, where 0=Prime, 1=2° prime and so on, an O(x), O(x$^3$) efficient algorithm is attached describing the process of factor generation and primality testing. The factors of each element are also available as the linear binomials are linked to the number of cycles per mod 18 in the non-prime matrices as row = a+1, column=b+1, in ((LOPI$_1$ + 18(a))((LOPI$_2$ + 18(b)), a,b≥0. The residues 3,15,9 are not coprime to mod 18 and are not needed as they do not generate a mix of prime and non-prime elements partition into non-prime matrices, all elements are non-prime that could stem from the non-coprime residues 3,15,9 mod 18. Once introduced, numerical analysis will be done with the system on a variety of long held definitions and conjectures using the continuity of the core identity and Dirichlet's Theorem on the infinite elements of primes in linear progressions, mq + b, where b and mod m are coprime, as in our system, the least odd partition identity/reduced residue system mod 18, with variable a=q≥0 and the 6 reduced residue as constants, b, yields the odd natural number set capable of containing the primes (LOPI($_{1,5,7,11,13,17}$) + 18(a), excluding congruence class 3, the only odd prime capable of being prime not coprime to modulo 18. The reduced residue system modulo 18 combined with linear convolution, a system we call, Efficient LOPI Matrix Algorithm, ELMA-18, provides a potential solution to the P NP computational problem, providing a polynomial ask and a deterministic polynomial answer to the question: Is this integer a prime or non-prime and if non-prime, what are its non-trivial factors. All prime elements in each congruence class are partitioned out as they are the compliment to the total 3-4 non-prime subsets within each LOPI superset congruence class, revealing 6 patterns of primes within the 6 congruence classes.

**Introduction**

LOPI/LEPI Identity under mod 18 introduced in this work:

$\quad$ LOPI/LEPI  N($_{a,bc…n}$)= (a + b + c…n).,= N≡ $LOPI/LEPI\ (mod18)$.

The lowest odd partition identity is used as a core grouping value that carries through infinitely in each element of each congruence class mod 18, a value kept constant in the partitioned non-prime subsets. This system builds the Superset from Dirichlet Linear Equations, (LOPI + 18(a)) and then non-prime subset elements by linear convolution of these linear binomials, ($LOPI_1$ + 18(a))($LOPI_2$ + 18(b)), containing the same identity/residue throughout the convolution equivalence relation. 18 is a special modulus in that it is used as a LOPI identity additive element, as 0 +a=a, +18 doesn't change our grouping core identity but allows for the generating of the next integer in the progression of each congruence class and later continued factor generation. As well, 18* is an absorbing element for the LOPI in this modulus system, used in the linear equations by which to grow the value of the next element in the progression, while keeping the grouping integrity of the LOPI Dirichlet linear equation during the binomial convolution, allowing us to generate all elements in each congruence superset and the factors of the subsequent composite subsets while maintaining the combinatorial relationships. The complement to the complete composite total subsets is the Prime subset in the Superset since, upon the convolution of the linear binomials, the primes cannot be mapped by the equivalence relation, as no multiplication by 1 is ever needed or allowed. The 3-4 images in the natural number codomain of the convolution for subset generation are the non-primes, the "anti-images" are the prime elements in each superset. The primes will not be represented in the subsets and upon query of the subset data matrices within each LOPI congruence classes, a 0 will be the result, only for primes will this be true. They cannot be represented and yet all non-primes within the superset must be generated by the sequential factor generating binomial linear equations that form the products in the subset data matrix. This is the basis of the attached algorithm for a new deterministic prime testing algorithm based upon linear binomial convolution whose polynomials are degree one. The resulting polynomials yield the composite elements of the subsets and the gaps in the codomain of the natural number range of the convolution are the prime numbers. Degree of primality, of a non-prime integer can therefore, also, be defined by multiplicity within the 3-4 proper subsets of each reduced LOPI, with 0 presence equaling Prime, 1=secondary prime, 3=tertiary prime and so forth. Degree of Primality equals presence in the 3-4 congruent LOPI subsets and the number of factors equals

Multiplicity(2) + 2, where +2 is for the trivial, 1 and itself. For labeling of a secondary prime we see the integer has only two unique prime factors as generated in the data function table. A multiplicity of one for secondary primes is always true by this definition in our data matrix, the only exception is LOPI$^3$, which is obvious and at the beginning. A LOPI cubed will result in a secondary "prime" in the sense that it is a unique case where the prime LOPI is multiplied by LOPI$^2$, thereby making the cube appear as a secondary prime since all the LOPI are prime, (except 1, which is always one no matter the exponent). For example,
((7 + 18(0)(7 + 18(2))=7X49= a LOPI class 13 element appearing one time, 7X7X7=343, 343 upon query yields a find of 1, indicating a multiplicity of one, but in this one case, $x^3$, the multiplicity is correct, it does have only two factors, but the idea that upon unique factorization secondary primes yield two unique factors is incorrect. The definition of a secondary prime therefore is an integer that has two unique factors that are not divisible by each other. The LOPI Prime cubed is not a problem for the system here as we are primarily focused on the partitioning and identification of the prime number subset within each LOPI class and then large secondary prime number factor generation. The unique LOPI$^3$ secondary prime is still a non-prime and will never yield a zero upon query of the subsets, and if we are searching to verify large secondary primes, LOPI$^3$ or LOPI$^3$ does not come into play. Any product from the factors generated producing elements in the LOPI non-prime subsets in this sequence, LOPI$^{3+1(a)}$, will never generate a multiplicity of one and never be confused as a large secondary prime. The multiplicity of such a number, LOPI$^{3+1(a)}$, will be ≥2 as LOPI$^{3+1(a)}$=LOPI$^{x-1}$(LOPI), exp. x=(3+1(a)), a≥0. We note LOPI$^3$ just to be thorough on the definition of a secondary prime multiplicity equals one when an integer has only two non-trivial unique prime factors. The linear binomial results of a product equal to LOPI$^3$, are listed as unique false multiplicity one results, that are not the product of two unique primes but do have only two factors but are not a secondary prime, but are a secondary square prime, LOPI P X LOPI P$^2$. Therefore, we make a special note for the return of 1 by the prime testing/factorization algorithm as the multiplicity of these values, as they are a secondary prime in that they equal a LOPI$_1$ (LOPI$_1$$^2$) and therefore really do have only 2 factors, but are not secondary primes: 5^3=125, 7$^3$= 343, 11$^3$=1331, 13$^3$ =2197, 17$^3$ =4913. Larger values in this sequence, LOPI$^{3+1(a)}$ will never generate a multiplicity of one and never be confused as a



large secondary prime and all factors are generated for these larger elements as well. The multiplicity of such a number, $LOPI^{3+1(a)}$, will be ≥2 as $LOPI^{3+1(a)}=LOPI^{x-1}(LOPI)$, exp. x=(3+1(a)), a≥0 and all factors will be generated for these elements.

Other moduli coprime to its residues, such as 12, can be used to yield factor generation but without the overlapping identity of residue=to LOPI, the predictable relationships are more embedded and the core grouping identity is lost. (12(17) + 7)= 211, this is a prime element corresponding to 211 element in congruence class 13 mod 18, the result does not have the core identity in common with the residue 7 and therefore the system doesn't have a readily available extra partitioning equivalence and we can't make predictive basic mathematical templates from mod 12, as it is not an absorbing and identity element within the linear progression ( Residue S(12)+ 12(a). However if we use mod 18 for 211, we know that 2+1+1=4 and the equivalent residue mod 18 for odd numbers is LOPI 13 for odds, therefore we know 211 is a residue 13 element and that in (mq + u), u=13 and m=18 and therefore q=11, or (7+18(11)). Upon non-prime subset generation through partitioning by an equivalence relation involving the binomial linear convolution of specific LOPI Dirichlet residue progressions needed for LOPI 13, 211 does not appear and is therefore is a prime, with zero generation into the non-prime subsets. The functionality of this system is based on the constant relationship maintained between mod 18 and the LOPI residues, allowing us to view each integer as an unrestricted partition of its core LOPI, by linearly adding in steps increasing by 1, cycles of 18, LOPI + 18(a), a≤0, we cycle further and further away from the core while still maintaing the core with which to identify larger and larger elements in the congruence classes, LOPI residues= totatives coprime to modulo 18.

Combinatorial templates are constructed as well under the rules of Modulus 18 LOPI/LEPI addition to give the congruence relations needed to produce all pairs of possible elements in each congruence class, LOPI Prime + LOPI Prime=LEPI Even, or $((LOPI_1 + 18(a)) + (LOPI_2+ 18(b)))= LEPI + 18(c)$, as a congruence relationship: *LEPI +18(c)≡$LOPI_1$ + 18(a) (Mod $LOPI_2$ +18(b))*. This approach provides a system of equations that provide all the possible, LOPI P + LOPI P= LEPI N, Goldbach Partition solutions using the infinite congruence classes into which all prime odd integers exist, except 3, no solutions are based on residues 2 or 3 as prime numbers are not



generated in these classes, as they are not totatives of 18, they are not in the reduced residue system modulo 18. They are unique within their respective congruence classes in that they are the only integer in their congruence class modulo 18 to be considered prime. We call the combinatorial templates for the even integers, the Goldbach/Dirichlet LOPI mod18 See-Saw partitions:

*LEPI +18(c)≡(*LOPI$_1$ + 18(a)) ((Mod LOPI$_2$ +18(b)) with a$_i$=0, b$_{max}$=c, or b=c-1,

depending on c, the number of cycles of 18 iterated in LEPI N. Sophie Germain primes, 2p+1, and Mersenne primes, $2^p$-1, can be analyzed as well in the context of the reduced residue system mod 18 and the correlation to the core LOPI identity. The core is found to affect the possibilities of the LOPI Super-sets in which they will be found as elements. If we perform 2p + 1 on the LOPI residues, and knowing that the core changes will be the changes throughout the whole congruence class of elements, then we see the classes are limited for both groups of primes and that there are also two types of Mersenne primes that predict which type of perfect number they will make by $(2^{p-1})(2^p-1)$. Also we show Sophie Germain Prime integers and their corresponding safe primes exist in LOPI groups 5, 11 and 17 and no twin prime pair exists that are both Sophie Germain primes, since twin primes are LOPI sets [5,7]; LOPI [11,13]; and LOPI [1,17], and each set has only one LOPI that can be a Sophie Germain prime or its safe prime. As well no elements in LOPI 5, 11, 17 congruence classes are ever made by the squaring of a LOPI, therefore Sophie Germain Primes and resulting safe primes are all generated only by LOPI mod 18 residue classes 5, 11, 17 contain no elements with a perfect square root. Figure 2 shows the base seed operations for the linear binomials in the non-prime convolution from which all non-prime elements are generated, it can also be seen here that LOPI classes 1, 7, 13 are the only classes with square base seed operations.

We use the even distribution of prime numbers in reduced residue classes and the Prime Number Theorem conjectured by Gauss **[1**], proven by Poisson and Hadamard **[2]** and the Dirichlet Theorem **[3]** in conjunction to show the primes are present. Also, by Bertrand's Postulate, **[4]** a long chain of uninterrupted non-primes has been shown to **not** exist, since it has been proven that there is at least one prime between n<p<2n. The convolution of the Dirichlet LOPI linear progressions, generates all the non-prime elements in our system in



each LOPI class leaves the primes unaffected, experimentally we see that the primes are the gaps in the codomain of the range of the convolution equivalence relation.  In the two different LOPI linear binomials needed in the templates of the Goldbach/Dirichlet see-saw congruence relation, the states of  LOPI $N_1$, LOPI $N_2$ in the seesaw will have the four pair possibilities, NP/NP, P/P, NP/P, P/NP, and cannot exclude P/P existing for the length of each linear progression from a=0, to $b_{max}$=c, c-1, all four possibilities will exist, as we show later by the non-overlapping convolution results for the non-primes, which therefore also yield unique gap patterns in the range of the codomain, the primes.

 Using the factors in the basic mathematical operation in each seed lowest LOPI product we get the linear Dirichlet progressions from which to perform the convolution as each factor has 18 added to it to make new products that form the future non-prime elements of the specific LOPI congruence class.  Constant maintaining of the core LOPI during factor generation by linear binomial LOPI convolution gives us an equivalence relation that can generate unique pairs of factors and subsequent non-prime elements that contain the same core LOPI from within the larger LOPI class. The primes aren't formed in the convolution due to the fact that the subset non-prime product elements are all elements with factors >1, coprime to mod 18 and therefore multiplication by one is never allowed to occur, insuring no prime can map into the subset data matrices.  The primes result from not being equal to any product of the factors generated, which is expected as do not allow one as a factor in subset convolution. The values which require one as a factor, the primes, are partitioned alone in each specific LOPI congruence class after partitioning of the non-prime elements from the LOPI superset. This feature in the system allows that primes upon query of the non-prime subset data matrices will not be found and output for this query will be zero.

 Each LOPI linear Dirichlet progression creates a congruence class, 18(a) + LOPI,  each one originates from a unique LOPI origin and each consecutive LOPI element is generated by linear step wise increases of cycles of 18, the classes cannot, therefore, overlap exactly on their P/NP states, or else they would be in the same congruence class, originating from the same LOPI core, which is can not be the case.  The congruence classes are disjoint and originate from the unique LOPI totatives of mod 18 = {1, 5, 7, 11, 13, 17} + 18(a), with a≥1,



none of the classes will ever match perfectly, but do create non-primes through convolution with the same mechanism but with unique sets of LOPI Dirichlet linear progressions to populate each congruence class. When considered in this context, therefore,

The non-primes and the gaps in the convolution of each LOPI linear binomial are cyclic, as demonstrated in the non-prime generation in **figure 3**, for LOPI 13. As the linear progression that contains the primes and the non-primes increases toward infinity in a linear fashion by cycles of 18(a), the inhabited space for the non-prime "plane and solid integers", through linear binomial convolution, is growing at a faster rate, and therefore the corresponding gaps in this cycle are growing faster as well, which are the prime elements to be available for all LEPI elements, which are growing as well, in the linear progressions LEPI even residues mod 18, {0,2,4,6,8,10,12,14,16}. The primes will be available to create the larger even LEPI congruence class elements by addition as conjectured by Crisitian Goldbach, P + P=2N≥6, correlated to the number of cycles equal or less by one to the number of cycles used to generate the LEPI in the LOPI partition seesaw. For example, LEPI 2, 686, c=38, b=37 at least, we have 37 cycles of our linear progression by+18(a) in which the gaps were generated that represent primes. Therefore if we have two template pairs for LEPI 2, we have  LOPI$_1$ +18(0) up to LOPI$_2$ + 18(37) progression range in which to find our primes, which we know are there as the convolution products are cyclic by a greater distance than 1 for sure and these gaps represent the primes in the the linear progression, the gaps are there from the convolution to delineate the primes sitting in the superset for us to combine as we wish. for each pair of templates iif the distance grows greater with time between gaps and closer with time for non-primes, **figure 3**, for LOPI 13. As well there are always at least two seesaw congruence relations per LEPI congruence class mod 18, therefore increasing our supply of primes from which to meet the Goldbach partition requirements. The solutions to different criteria can be determined by the basic results of performing the desired mathematical operation on the LOPI residues, the resulting solutions are predictable relationships due to the core predictability. We know how the core LOPI will be affected, or how an integer will be manipulated, by a particular operation. The basic mathematical relationships between the LOPI/residues are very useful in this system for a variety of numerical inquiries. For instance, we will explore the LOPI



congruence classes that can contain Sophie Germain Primes and present a finer congruence relation that an integer must have to be a Sophie Germain prime. We also explore the Twin Prime Conjecture, as all possible twin primes exist as the non-prime subsets of LOPI 17.  Upon setting a=b for 5 +18(a), (7 +18(b) and for LOPI 11,13, and a=b+1, for LOPI 17,1, factors generated for convolution can be verified as resulting in an element with only one result, multiplicity of 1, in the specific 3 subsets data matrices of  the LOPI residue class 17, therefore indicating a secondary prime made from twin prime factors, since a and b were predetermined to yield factors separated by 2 in order to be twin elements in their respective LOPI. For example 323=19*17, or generated by convolution of the linear LOPI equations,  ((1 +18(1))(17+18(0))), a=b+1.  Another example, 57599=((5+18(13))(7+18(13))), or 239*241=57599, a secondary twin prime, located one time in LOPI residue subset 17, and with a=b=13, therefore made from twin prime factors, and are a twin prime pair.  Without partitioning integers based on the core identity of the LOPI, the totatives of mod 18, these basic predictable mathematical relationships would not be possible to make the combinatorial relationships that yield the information within this unique system, as 18 as modulo used also as an absorbing element in 18(a) and identity element by adding this neutral quantity to the core LOPI + 18(a). As for example ((7 +(18*17)), the core identity of 17 is absorbed by multiplication by 18,  to yield the neutral quantity of 306 in LEPI 18, and upon addition of the core LOPI,+ 7,  we produce again an element of LOPI residue class 7, or 7+306=313, odd LOPI family 7. We can now know within the congruence relations  what is possible and not and that the LOPI P + LOPI P  solutions will exist.  The Goldbach Partitions are shown to be possible as a consequence of the fact that primes exist where there are no factors to be multiplied by other than 1, to yield an integer of the value of the prime in the reduced residue classes, as the pattern of iteration by +18 to both, Left or Right and subsequent multiplication, has already created all the composite integers leaving sufficient gaps/primes to add to yield the even equivalence classes of modulus 18.  Iteration of 18 and subsequent multiplication cannot generate the integers that are created by multiplication by only 1, the primes, but does generate all the numbers by convolution in a congruence class that have odd factors co-prime to 18 >1, as the equivalence relation cannot partition out the primes contained in the Dirichlet linear equations of the Superset, while the factors for all composites are



generated by multiplication after iteration of +18 to the factors if the 3-4 lowest seed operations in each congruence class. We see an interdependent pattern, relationship, between the primes and non-primes. The generation of the non-primes, by applying an equivalence relation to the LOPI superset, leaves the primes untouched, as they are not created by multiplication by any other linear LOPI Dirichlet equation, they are made only by the original linear equation, that can be thought of as being multiplied by 1. Unique factor multiplicity is a useful tool, as well, in the sub-grouping of the non-primes into subsets. We show that an integer's degree of primality in the LOPI Supersets can be determined by the number of times it appears in the subsets. Moreover, it becomes clear that if an element appears more than once, with a multiplicity of >1, it has unique integer factoring that allows it to exist in different subsets, therefore primality can be defined by an integers multiplicity in their respective composite subsets. For example integer 175, is an element in the LOPI Superset 13, (17+18(9)), upon convolution from stem 85, 175=(5 +18(0))(17+18(0)), and is generated as well in subset stem 7X7, stem 49, ((7+18(0))(7+18(1))), we see a multiplicity of 2, which equals a tertiary prime, 2(2) + 2=6/2=3; 4 non-trivial factors and 2 trivial, self and 1. We define the degree of primality of the integer based on its generation into its core LOPI subsets, therefore 175 is generated two times by convolution in the subset or subsets of Residue Class 13, therefore is a tertiary prime. The number of factors equals number of times found in subsets times 2 and plus two for trivial factors since a non-prime generates into subsets only by unique factor multiplication, 2(2) + Trivial 2=6, 6/2=3 Degree of Primality. The fact that these integers exist in different places allows us to treat them uniquely in their own space, which allows for primality labeling based on generation into well defined sequential non-prime subsets. A new definition of a prime or a non-prime could be introduced, as a prime being an integer that exists only in one linear space, as defined by its one generating linear Dirichlet Equation under the reduced residue system modulo 18, while non-primes, or composites, can be defined as integers that exist, as Euclid reported of Nicomachus **[5]**, "a prime number..,according to Nicomachus, because it can only be arrived at by putting together a certain number of units," a collection of one units in linear space, or as shown by r.r. system mod 18, the Dirichlet Progressions, LOPI + 18(a), contain all the prime odd elements possible, excluding the one 3, and that upon applying an



equivalence relation to generate the non-primes, the primes are left untouched. They remain the only subset within each LOPI superset that do not map into the proper non-prime subsets. By partitioning the natural odd numbers into the 6 congruence classes coprime to mod 18, which coincides with the core LOPI identity, the linear eqaution that produces the whole LOPI Superset and upon subsequent linear binomial convolution to produce the non-primes, the primes show to be immune to the equivalence relation of convolution. The convolution that captures all the non-primes within each LOPI congruence class and places them into 3-4 subsets, captures their multiplicity based on unique factor generation and in so doing we see the contrast between the prime elements and the composite, as reported by Euclid, again, expressed by Aristotle,[5], "who contrasts the composite number with that which is only in one dimension." We find the prime elements to be only in the LOPI reduced residue supersets, which are each generated only by one of six unique linear Dirichlet progressions using the reduced residues, LOPI Mod 18, "a prime number is that which is *linear*"[5],. The two variable linear polynomial that represents the convolution leads to the generation of only the non-prime elements and their factors, the primes will have no inverse in the function. The same holds true if the element is later generated again within the same subset, still we retrieve its presence in the subset and count its multiplicity, since its generation will be due to unique pair factors, this multiplicity indicates the element existing in another unique subspace, or a later point in the same subspace. More over the location of the element in the data function table is linked by equation to its factors, where a and b in the linear equations are row-1 and column -1 in the linear convolution $((LOPI_1 + 18(a))(LOPI_2 + 18(b)))$, therefore upon query we can note the row and column, indicating the factors and generating those specifically if needed.

Discussion

The reduced residue system mod 18, introduced in this paper and referred to from here out, as the LOPI mod 18, is a core quantity used to partition all of the odd natural numbers that can contain primes, as well as non-primes into the reduced residue system modulo 18, the lowest odd partition identity, the reduced residues of mod 18, corresponds as well to the totatives of ω(18), reduced residue system modulo s=ω(m) **[6].** Upon



rearrangement of the number line 0-18, **figure 1,** the lowest odd partition integer is seen to be a core constant grouping identity for the partition of the odd natural numbers coprime to mod 18, forming the six LOPI congruence classes. Each LOPI has a corresponding LEPI as well, but to partition the odd natural numbers into 6 congruence classes coprime to 18, we will use only the odd residues coprime to 18. The LOPI is a unique and constant identity in each odd number under modulus 18, defined by $N_{LOPI}$≡LOPI (mod 18), as well as equaling the sum of all the digits of the partition integer, a+b+c in 311, equals LOPI 5, and 5 equals LOPI in 311≡5(mod 18). Each set of a LOPI/LEPI family is populated by iteration of +18, much as addition by 0 doesn't change a, a + 0 = a, LOPIn + 18 = LOPIn', the LOPI grouping identity doesn't change, therefore the system works due to the use of 18 as a type of additive and absorbing element. The Dirichlet seed progressions use this quality also, in that 18 is equivalent in function to the absorbing element 0 in the number system, 0(a)=0, 18(a)=LEPI 18, no matter the value of *a* in the linear equations, we are just adding a neutral LEPI 18 quantity. Essentially the linear equations become extensions of the core grouping LOPI and populate each of the 6 congruence classes, while maintaining the integrity of the grouping core LOPI. LOPI + 18(a)=$LOPI_{n1}$, the progression continues with *a* being the set of natural numbers, a={0,1,2,3,…}. LOPI set= {LOPI + 18(0), LOPI + 18(1), LOPI + 18(2),…}. All odd non-prime natural numbers and their factors are generated as well from the convolution of a set of Dirichlet binomial equations that stem as well from the totatives of 18, or the reduced residues of modulo 18; 1, 5, 7, 11, 13, 17.

The LOPI/LEPI identity, corresponds to congruence identity relationship of the reduced residue system of Modulo 18 as:

LOPI/LEPI Identity under mod 18:

$$\text{LOPI/LEPI } N_{(a,bc..n)} = (a + b + c…n)., = N \equiv LOPI/LEPI \ (mod\ 18).$$



Where LOPI/LEPI= (a + b + c...n)., denotes summing down to lowest odd partition integer which will be one the integers less than 18 in the complete covering system: 1, 3, 5, 7, 9 11, 13, 15, 17 for odd integers, or LEPI for even integers: 0, 2, 4, 6, 8, 10, 12, 14, 16, and mod 18(a). Modulo 18 makes a complete disjoint irredundant covering system **[7]** for each congruence class upon iteration + 18 of each element and each mod by 18(a), since the grouping identity is based not only on the modulo but also on the one lowest odd/even partition identity contained in all integers. This initial modular grouping identity, that continues every +18, or LOPI + 18(a). This quantity of 18 is found to be key in factor generation that is later used to populate the non-prime elements within each of these congruence classes on the reduced residue system modulo 18, the identity uses mod 18 since we keep mod 18 need constant in the equations, as it always returns the original least odd or even partition identity for every natural number on par with the parity of the integer. Division by 18 yields the even LEPI for an even integer and yields the odd LOPI for the odd integers in the reduced residue classes of mod 18, therefore the other modulo of 18(a) are not used in the identity definition of the LOPI/LEPI but in the linear binomials, 18(a) is an essential feature as it provides a neutral quantity through multiplication of 18(a) to add to the LOPI that will not change the partitioning core LOPI identity, LOPI **+ 18(a),** a neutral quantity (LEPI 18) that acts as adding 0 to the LOPI constant, allowing for the next element in the congruence class. Each integer is the LOPI plus consecutive cycles of 18, therefore each element in each LOPI congruence class meets the congruence relation $N_{LOPI} \equiv LOPI(mod\ 18)$. The other elements of 9/18 can not be used to consistently return the LOPI in the reduced form of the LOPI and LEPI, as they are not all coprime to all LOPI mod 18 and therefore cannot be reliable as a basis from which to create combinatorial relationships.

The congruence relation LOPI N≡LOPI/LEPI Mod 18; {0,1,2,3,4,5,6,7,8,9,10,11,12,13,14,15,16,17}, {$a_1$(mod $n_k$)} , {$a_i + n_i x : x \in Z$}, due to the identity LOPI/LEPI= LOPI/LEPI N≡LOPI (mod 18), the whole number line, [0] [1] [2] [3] [4] [5] [6] [7] [8] [9] [10] [11] [12] [13] [14] [15] [16] [17], can be represented in the form of the LOPI/LEPI progression for the odds and evens: LOPI + 18(a), and the linear equation: LEPI + 18(a), for the even integers. Set Z is understood to be all LOPI mod 18 classes containing odd primes and non-primes. Any number, large or small, can immediately be placed in one of the 18 congruence classes with minimal effort, giving us valuable



information and when we used the reduced residue system that can contain all odd prime/nonprime integers coprime to 18, we know where to look for prime information and factor information from the data function tables produced upon generation of the non-prime subsets within each LOPI superset.

The premise that all odd/even integers have a constant core LOPI identity, is evidenced by the modular arithmetic of modulo 18, shown in **figure** 1. Upon addition of 18 to the reduced residue system of mod 18 and as well to the modulo 18 we see that each odd element in the LOPI congruence classes contains the LOPI as well, as we are just adding cycles of 18 to a constant core quantity. This core quantity is lowest odd partition identity that corresponds as well to the totatives of ω(18), to form the reduced residue system mod 18, introduced in this paper and referred to from here out, as the LOPI mod 18. All the odd natural numbers that can contain primes, as well as non-primes are contained in this reduced residue system. The lowest odd/even partition identity is one with infinite elements grouped into these six congruence classes, as each positive integer has only one LOPI/LEPI as a core grouping identity. The LOPI congruence class sets of {1, 5, 7, 11, 13, 17}, for modulo 18 correspond to the 6 odd totatives of *ω(18)=6*, the Euler function for the set of numbers coprime to the modulo, as shown in Theorem 4-5, "if s integers $r_1, r_2,....,r_s$ form a reduced residue system, modulo s=ω(m). **[6]**. By viewing each integer in each LOPI superset as a larger version of its kernel LOPI we create a relationship between these groups that can be utilized to generate the factors that lead to the composite numbers in each class, while still retaining the original grouping LOPI identity. From here on we will only be dealing with the reduced residue system of Mod 18 as they are the only congruence classes forming the basis for sub-grouping all composite elements under each congruence class coprime to mod 18, these will be known as the LOPI (mod 18).

To generate every element within the 6 LOPI Supersets, we use the Dirichlet linear equation in the form of *LOPI + 18(m),* or *N = m(q) + LOPI*. The value of each LOPI is relatively prime to 'm', 18 and is used as a constant in the linear equation for the superset and in the structure of the linear binomial convolution that generates the elements of each non-prime LOPI subset. Given the definition of a Dirichlet Progression**,[8]** there are an infinite number of elements of primes generated by the linear equation in the infinite superset, LOPI + 18(a)



and given the Prime Number Theorem, **[9, 10]** it has been proven that the these prime elements are distributed evenly within the Dirichlet arithmetic progressions where g.c.d. (u,m)=1, such as we generate here for each congruence class with the reduced residues of mod 18.  Equivalently, N LOPI ≡LOPI *(mod 18)*, or by Euclid's division lemma [13], *N = m(q) + u, N LOPI=18(q) + LOPI,* **[11]**. In relation to Theorem 4-4 **[6],** where u=LOPI, m=18, q=coefficient, and the LOPIs=ɷ(18)= 1,5,7,11,13,17, form the Complete Reduced Residue System, $N_{LOPI}$≡*LOPI(mod 18).*  The group order for cyclic group 18 **[12, 13]**, a 2p^2, is 6, with the generators, g^k being 5 and 11 for the all residues of mod 18, not all the elements in each LOPI, which repeats every 6(a). LOPI congruence classes 5, 11 identity elements, e, corresponds to $g^0$ the Euler totative=to the LOPI=residue after division by 18.  The identity element residue cores of LOPI $C_6$=<g>=the generator, $g^{0+6(a)}$, $g^{1+6(a)}$, $g^{2+6(a)}$, $g^{3+6(a)}$, $g^{4+6(a)}$ $g^{5+6(a)}$ return the residue elements of g=5,11 as capable of generating all 6 residues mod 18 in all LOPI congruence classes, the finite cyclic group of the LOPI reduced residue system mod 18, {1,5,7,11,13,17}. Therefore residue 13 for instance can produce only non-prime elements by exponentiation in LOPI 1,7,13, or the non-prime Non-Sophie Germain Residue classes, $N_{LOPI}$≡1(mod 6) and residue 1 returns itself no matter the exponent, 17 returns only elements within non-prime subsets LOPI 17 and LOPI 1 elements and 13 also generates non-prime elements in the non-Sophie Germain LOPI classes, 1,7,13, as does residue 7, only non-prime elements in LOPI classes 1, 13, 17.  Therefore we can say that LOPI non-prime elements in LOPI  5, 11, 17, possible Sophie Germain Prime classes, $N_{LOPI}$≡5(mod 6), are never non-prime LOPI [2] elements, they are never squares, they have no square roots as none of their non-prime elements are generated by any of the residues LOPI 18 squared.  They are generated only by $5^{3+6(a)}$, $11^{3+6(a)}$, $13^{1+6(a)}$, $17^{1+6(a)}$. The Euler totient function of 18 is 6, as expected in modulo 18 with 6 reduced residues, group order is 6, two generators capable of returning elements within these 6 classes at cycles of 6(a), g=5, 11. We are producing the 6 reduced residue classes of mod 18 where the residues equal also the core LOPI identity, which repeats every 18 from the core LOPI totatives of 18.   The 6 LOPI reduced residue congruence classes mod 18 contain all natural odd prime integers, as well as all non-prime integers coprime to mod 18, this is our Universal set of integers of interest.  LOPI 3



residue class mod 18 is the only element considered prime not coprime to 18. It doesn't interfere with the system as beyond 3, $N_{LOPI} \equiv 3 \pmod{18}$ generate only non-prime integers.

In the LOPI mod 18 system LOPI is used to create 6 super groups of the universal primitive group, of the Universal set by creating 6 separate disjoint Supersets based on the residue identity of the totatives of modulus 18 with the linear Dirichlet equation *LOPI + 18(a)*. The LOPI of a number, which is the constant residue identity also determines the linear Dirichlet equations in the convolution to generate the factors of the elements that populate the composite subgroups of each reduced class.

In keeping the residue mod 18 a constant residue grouping identity, which also corresponds to the lowest partition identity of odd natural numbers, we use these overlapping constants as a basis for element generation, both prime and composite. For each residue u, there is an integer N mod m, such that LOPIN≡LOPI(modm), that allows NLOPI to also have the LOPI as a residue identity. Both are the same value of the totatives of ϖ(18) by which they were originally grouped into the 6 separate supersets, LOPI + (18a), with a≥0. In using the reduced residue mod 18 as a constant by which to group all the odd natural numbers, with r always being coprime to mod 18, or the totatives of 18, each LOPI set and its subsets are complete and congruent under mod 18. No grouping is based on LOPI 3,15,9 or on the even residues of modulus 18, therefore we have formed the primitive multiplicative group of odd natural integers coprime to modulo 18, grouped them into six infinite multiplicative groups with corresponding non-prime subset groups, all coprime to modulus 18 and congruent within their LOPI family, with the congruence relation, (*a= nk + b), 0≤b<n* and g.c.d.(n,b)=1, where b=LOPI residue (mod18), therefore *b* equals the LOPI remainder after the division *a* by 18. Each element in the natural odd number set coprime to 18 exists uniquely in one LOPI superset as there can be only one constant grouping residue of modulo 18 for the elements in the LOPI congruence class super-sets and in their respective subsets, even in a non-reduced systems, division of a number by another less than itself dictates that only one remainder, in general, will be derived from division. This isn't to say that the elements in the subsets will only have one pair of unique factors generated by the convolution, but it is to say that each non-prime element in their respective LOPI subsets will have only one LOPI reduced residue mod 18, therefore,



All LOPI supersets/subsets are disjoint within each congruence class to the other LOPI Supersets/subsets, as the subset elements are grouped under the same LOPI identity residue as its parent LOPI superset via mod 18.

By partitioning based on the reduced LOPI identity, all the elements of each superset have a constant that can be used to predict the effect of some mathematical operations on elements in the same LOPI congruence class. The effect on the residue constant integer is seen to be the same effect on any larger LOPI integer, carrying consistently throughout the elements of the same congruence class. The LOPI/LEPI modulo 18 system lends a structure to the set of natural odd number Dirichlet progressions, LOPI + 18(a), and upon convolution of specific linear equations, reveals the pattern of factor generation of non-primes, where one as a factor is never allowed, revealing why the primes remain in the superset upon convolution, the primes' only factor besides itself, one, is not allowed in the convolution and therefore the primes are not produced by convolution in (LOPI$_1$ + 18(a)(LOPI$_2$ + 18(b), with a,b≥0 in the general case. **Figure 2** shows a Venn diagram of the seed factors that are used to generate the 3-4 proper subsets of each LOPI superset, which contain all odd non-prime natural numbers of that LOPI. The intersection within LOPI subsets is not shown, but we note here that overlap does occur accounting as well toward the unique factor pair product multiplicity upon query of each subset data matrix, discussed later in the presentation of primality testing by unique multiplicity, an algorithm based on non-prime subset matrix query for obtaining factors as well as determining prime degree of a natural odd number in efficient $O(x)$ to $O(x^3)$. Factors for natural odd numbers are coprime to 18, i.e., not divisible by 3. The complement of the subsets in their superset are the primes within that LOPI. We define 3-4 equivalence relations in order to partition out all the non-prime elements within each of the 6 LOPI reduced residue super-sets, with the last subset remaining in the super-set, being that of the primes: By theorem 2.3.4, **[14],** which reads:

"Let R be an equivalence relation on the non-empty set LOPI A. To each *a* member A define â to be the subset $\{x: x \in A \wedge xRa\}$ of A. [Thus â is the subset of A comprising all elements x of A which are related to *a* under R.] Then the subset $\zeta$ is a partition of A."



The prime subset within the mixed Prime/non-prime LOPI superset, A, will be partitioned out by the subtraction of the set of composite subsets, ζ=non-empty set of non-prime subsets, whose elements, x of A, xRa, are generated by the linear Dirichlet binomial convolutions for each LOPI congruence class. The prime set is unique and not in the elements partitioned into ζ= {â : *a ε A}*, where â is the subset of all non-prime elements x in the LOPI super-set, A, related to *a* under the e.r. R, the convolution of the LOPI specific linear binomials, to yield the set of NP subsets, {*x* : x ε *A* and *xRa*} of A. Each subset has its own equivalence relation which partitions the 3-4 non-prime subsets in the complete ζ= {â : *a ε A}*. LOPI Superset\LOPI Subsets = LOPI Sub' = LOPI Prime subset within A, the reduced residue LOPI superset.

The complement of the subsets in their superset are the primes within that LOPI. We define 3-4 equivalence relations in order to partition out all the non-prime elements within each of the 6 LOPI reduced residue super-sets, with the last subset remaining in the Super-set, being that of the primes: The primes do not undergo an equivalence relation, R in A, and are therefore remain the only subset unchanged in LOPI Superset A, while all composite elements have undergone a transformation by binomial linear convolution and are partitioned out of A under R into three non-prime subsets for LOPIs 5, 11,17 and four subsets (2 whole and 2 half sets, equal to 3) for LOPIs 1,7,13. Interestingly the LOPI classes that contain four subsets of non-prime elements produced by a LOPI * LOPI, or LOPI$^2$, have four subgroups and correspond to the LOPI that cannot generate a Sophie Germain prime, LOPIs; 1, 7, 13, as will be discussed below.

At this point we have 6 infinite disjoint super-sets containing all odd natural numbers of interest, possible prime and difficult to factor composite numbers grouped under a common identity, LOPI, and from here we can derive the composites as we go by generating the factors and subsequent product of all the composite numbers within each LOPI superset. $A_{super}=B_{sub}+C_{sub}+D_{sub}+P_{sub}$, or A- [subsets $\zeta = $ â: $a\epsilon A$} ]= Subset P, the primes are the complement of the subsets in the superset. We create an infinite superset with linear binomial Dirichlet progression using the reduced residue system mod 18 and further partition the Superset into subsets by the convolution of reduced residue linear binomials mod 18, whose equations within the convolution keeps the core LOPI superset residue identity intact for each product element of the subset. No element is the result



of multiplication by one, as the system restricts this operation. Each LOPI composite subset will be defined below in the more detailed discussion of set and subset generation. The infinite absolute complement is the infinite set of elements, not in U, or *P set=U\Total LOPI sub,* where U= Universal set, containing all the elements of the 6 LOPI mod 18 sets of natural odd numbers. This infinite absolute complement, Total {*LOPIsub'}={P}*, the infinite set of all odd prime elements:

{$U_{r.r.mod\ 18}$}-[Total Subsets $\Omega$ ]= [$\zeta_{all}$= {â : *a ε A*}] = {P}.

**LOPI reduced residue mod 18 Set and Subset Generation**

We reduce the natural numbers, N, as a partition of the lowest even or odd partition integer, using modulo 18, seeing the infinite number of ways in which the lowest partition integers can be formed, we deterministically find the factors of an odd integer coprime to 18 and assign unique multiplicity, which is a measure of primality in the non-prime LOPI subsets. By using the unrestricted partition definition of an integer, each LOPI N element (partition of the LOPI) is in its LOPI/LEPI linear progression mod 18, and as well for the non-prime, being generated by the linear binomial convolution function that originates from one of the LOPI seed non-primes. These elements become "planar" or "solid" elements, as discussed in the Elements book VII [] by Euclid, prime numbers can be viewed as only linear, having no length or width, while composites numbers are divided into plane numbers and solid numbers.  By not putting a partition restrictions on the partition definition of number of parts, we allow that each LOPI/LEPI congruence class has an infinite number of elements with the LOPI core identity of one of the totatives of 18, LOPI N=($LOPI_A$ + 18(a)) as *a* goes from 0 to infinity.   These unique initial non-primes with factors of the totatives, $\omega(18)$, or the LOPI, are the basis from which all other composite integers generate, or elements in the LOPI NonPrime Subsets, "grow".  All the composites within each LOPI congruence class can be generated from the 3-4 lowest non-prime seed elements in each LOPI Congruence Class Superset and upon division of mod 18 the LOPI of that composite will always equals the reduced residue identity of the LOPI congruence class.    LOPI/LEPI  $N_{(a,bc..n)}$= (a + b + c...n).,= N≡LOPI (mod 18), for summation of



digits, go down to lowest residue mod 18, even to even, odd to odd. For example 1111 is an odd number in LOPI class 13, therefore although it sums directly to 4, we must use 13, as that is the residue mod 18. On the other hand, 74 sums to 11, but it's an even LEPI 2 element, therefore it has 2 as its core LEPI residue. Each element of the residue classes have an even version and an odd version, making up each 9 residue classes that represent our seemingly symmetric natural number system, 0/9, 2/11, 4/13, 6/15, 8/17, 10/1, 12/3, 14/5, 16/7. Their coexistence and constant identity in each number leads to many mathematical combinatorial controls that allow us to create specific types of numbers. Here we will show how we used the LOPI to generate all the non-primes within each LOPI superset in order that we may partition them off, to let the primes show in the LOPI superset.

In the number line figure 1, the seeds are the lowest non-prime elements made from the reduced residue of mod 18 in each progression. In figure 3, we show the iteration by 18 through for LOPI 13. In the ternary tree the duplicates begin to occur as shown in gray. During linear binomial convolution none of this occurs as it is point-wise column vector times row vector multiplication and redundancy is avoided, allowing for the unique multiplicity in the system mod 18.

The proper subsets are each generated by the convolution of two linear binomials, derived from the original LOPI linear Dirichlet equation. The two linear binomials are the single ones that make up each LOPI linear equation, $((LOPI_1 + 18(a)(LOPI_2 + 18(b)))$, the seed's factors determine which two we use to generate the non-prime subsets. After that in the convolution we just iterates 18 every time and we multiply pointwise each new factor by each old one.

Commonly, consecutively adding 18 to each LOPI populates each column of Supersets, N=LOPI+18(a), a≥0, as seen in the diagram 1. Each LOPI sequence is an arithmetic infinite progression from the Dirichlet form u + kd, where u=LOPI and d=18 and k≥0, to yield the superset:

*Nlopi={lopi + 18(a), lopi+18(a+1), lopi +18(k+1),….}, with $X_n$ = LOPI + d(n-1),*

For the purposes here, it is better to rewrite $X_n$ in terms of $X_a$, by the relationship derived between *a* and n; a=n-1, n=a+1**.** Putting these equations in terms of the coefficient of 18 facilitates Row and Column value ,as 0



is a coefficient for 18 in the convolutions, except for LOPI 1. The LOPI value is included in the Superset where a=0 in $X_a$.

Therefore;

$$X_a = LOPI + 18(a) \text{ and } {}_{x0}\sum{}^{xa} = (a+1)/2(2(LOPI)+18(a)), \text{ from } LOPIa_0\text{-}LOPI_{aN}$$

For example: N=125= 125/18= (17 + 18(6))=125= 6 residue 17 for $X_a$=125 and ${}_{x0}\sum{}^{xa}$ from $LOPIa_0..LOPI_{aN}$ =

(6+1)/2*(2(17) + (18(6))= 497 and ${}_{x0}\sum{}^{Xa}$ $X_{0,...}X_6$ = 17+35+53+71+89+107+125=497

All the properties that apply to a linear congruent infinite progression under modulo 18 applies to each LOPI Progression with each new line of +18 creating a complete reduced residue system via mod 18(a); the new Modulo of that cycle will return back the same reduced residue system as well, **[15].**

LOPI={LOPI+18(a), LOPI +18(k+1…)}, where a≥0, which is in our linear binomials the modulo system 18(a).

While this is used for the linear progressions and the linear non-prime generation, the LOPI 18 must be used in the LOPI identity equation, as 18 always returns our original LOPI residues, as they are the totatives of 18 and not all later moduli 18(a) have the same residues as mod 18.

Each LOPI Superset, populated by $N_{LOPI}$ =u + k(d), with Integer LOPI N=LOPI + 18(a), a>=0, and its proper subsets are disjoint to all other LOPIs and contain no elements from any of the other LOPI supersets or their subsets. Integers can have only one LOPI/LEPI core, as it is derived from the identity: $N_{LOPI/LEPI}$≡LOPI/LEPI (mod 18)=Sum a+b+c..of Integer N and all subsets elements in the same LOPI Super-sets carry forward the same core identity upon convolution.

$N_{LOPI}$ Subset=LOPI(mod18), generated by ($LOPI_1$ + 18(a)($LOPI_2$ + 18(b) ), a,b ≥0.

The uniqueness of mod 18 is the fact that all members of the superset and subsets have the same LOPI upon division by 18.

18 is an absorbing and identity element in the LOPI mod 18 congruence system, $LOPI_0$+ 18(a) = LOPI $N_{a+1}$ and 18(a)=LEPI 18, an absorbing element, therefore the combination yields a neural quantity that leaves the core grouping identity in tact for the elements generated by the linear Dirichlet equation and later by the



convolution, as we multiply two LOPI linear equations to yield the desired LOPI non-prime subset elements, ($LOPI_a$+18(a))($LOPI_b$ + 18(b)) yields LOPI $N_{a+1}$*LOPI $N_{b+1}$, a non-prime LOPI element for the congruence class subset belonging to the same congruence class as the superset. LOPI 13, for example, superset elements are generated by {(13 + 18(a), 13+18(a+1),..} and the composite subset from the convolution of the seed linear LOPI equations that yield 85 in LOPI 13, 17*5, is ((17 +18(a)(5+18(b)). The point-wise convolution of this equation with a,b≥0 gives all non-primes that can be generated from this seed in the LOPI superset 13.

Subset generation:

Each LOPI's subsets elements are mutually congruent under Theorem 5-1 in Andrews, G.**[16]**
If d=g.c.d.(a,c), then the congruence an≡ b (mod c) has no solution if d doesn't divide b, and it has d mutually incongruent solutions if d divides b. Since we have a g.c.d. of 1 between all *a* and mod 18, then all solutions are mutually congruent in the ($LOPI_1$ + 18(a))($LOPI_2$+18(b))≡LOPI(mod 18) for each LOPI residue.
 Following Andrews example, the subsets in the form of their seed operations give solutions that are all congruent as Andrews gives an example of g.c.d. 1, since we designed our system to create only the LOPI elements that are coprime to modulo 18, we can therefore present the case that all solutions in each LOPI subset are are mutually congruent.

The first non-prime seed operation generates "node" factors for LOPIs of that subset. For example the LOPI 13; (17)(5)=17(mod18), g.c.d (17,18)=1, 85-17=72/18=4, or 85=13 +18(4), (17 + 18(0)(5 + 18(0)),
 this becomes the binomial generating function for each element in this LOPI 13 subset originating from 85, where a≥0, b≥0. By using the LOPI itself and 1, with b ≠ 0,and the lowest first unique non-primes as seed operations in each LOPI, all non-prime elements within each subset of each LOPI are generated to reveal a data function table showing each non-prime and its factors **(Appendix A)**. Figure 3, below shows this growth process with a ternary tree of LOPI 13.



All the sets generated are proper subsets of the superset A, and the Superset is a subset of the

{Universal set}= [{LOPI 1}, {LOPI 5}, {LOPI 7}, LOPI 11, LOPI 13, LOPI 17}].  **figure 2**

As with identities certain attributes can be attributed to each LOPI grouping. All elements within each LOPI superset are congruent and the LOPI + 18(a) element will be represented only once in each superset, no overlap between LOPI congruence classes is possible, as they are populated by +18 and each number has only one lowest odd partition integer, as found by the remainder of division by 18 or $N_{LOPI} \equiv u \pmod{18}$. The LOPI sets are disjoint from the other LOPI sets and each one can be viewed as separate systems by which to analyze a specific integer based on its LOPI. The binomial expressions used to generate each non-prime within each LOPI are specific in that a LOPI can only be calculated in specific ways using the other five LOPIs, or the lowest non-prime integers in each set are made of other LOPI that when multiplied maintain the LOPI of the congruence class. For example, the generating seed operations form the binomials that form the non-primes within each group, LOPI 13 by multiplication can be formed by 4 equations, therefore will have four proper subsets consisting of only non-primes, as there is no multiplication by 1, since b≥ 1 in 18(b) for the non-prime generating arm of LOPI*1. All solutions have factors greater than 1 and therefore have no possibility of being prime by the definition of prime as a natural number with factors of only 1 and itself.

As an example of the linear seeds used to generate the NP subsets in LOPI 13:

13x1, 7x7=49, 11x11=121,  17x5=85.

All these seed operations consisting of the other LOPI and reveal that each non-prime grows from these first parent reduced residue (r.r.) products and are the lowest non-primes generated in the linear progression LOPI + 18(a), a≥0.  In congruence class LOPI 13 each non-prime generated from the seed operation of 17x5 fits the binomial expression (17+18(a))(5+18(b), with a≥0,b≥0.  As shown in **figure 3,** the seeds of LOPI 13:

   13 x1          7x7          11x11       17x5



Each of these four seed operations multiply to LOPI 13 and upon adding 18 to either one or both and then performing the binomial convolution, the LOPI is unchanged due to the fact that addition in this Mod 18 LOPI system is like addition by 0 and multiplication by 18 to any coefficient of 18 is equivalent to multiplication by 0 in the natural numbers with the residues in mod 18, therefore we have essentially addition by 0 when we have +18, in that the LOPI is unchanged. Multiplication of a LOPI r.r. by 18 absorbs the residue grouping identity and leaves the equation with addition by "0" or 18(a) + LOPI, for the purposes of the grouping LOPI identity. The value of the number changes, but the constant residue identity mod(18) stays the same. The beauty of multiplication by 18 as opposed to 0 gives the ability to absorb the identity of the coefficient a, in LOPI + 18(a), without destroying the value of the integer to 0, leaving a quantity that can then be added to a LOPI that doesn't change the LOPI identity, since we are adding by a multiple of 18, but leaves the next larger element in the LOPI congruence class.

Additive identity element within the linear Dirichlet equations, LOPI (+ 18(a));

 LOPI/LEPI system mod 18, $N_{LOPi} \equiv$ r.r.LOPI (mod 18) , +18(a) doesn't change the grouping LOPI identity partition.
Multiplicative element LOPI/LEPI system mod 18,; 18*(a), with a ≥0.
LOPI + 18(a), 18*a, always equals a neutral even quantity, since even 18(a)=even integer, we then see that;
 LOPI odd + even18(a)=LOPI N , or, always a new odd member of the LOPI congruence class.

For Example, LOPI 13 non-primes are generated when 18 is added to either one of the multiplicands, the L, both or the R of the seed operation, then the convolution equation is reached, in the form of both, with different coefficients for 18 in each linear equation being convoluted. The use of addition and multiplication by 18 is demonstrated below allows the elements to grow without changing their grouping identity, either in the Supersets or the consequent subsets resulting from the linear convolution.

  17 x 5 = 85, each factor is one of the  r.r. of 18 and produces one of the lowest NP elements of LOPI 13, under the linear convolution equation:  ((17 + 18(a))((5+ 18(b)), where a≥0,b≥0.



|     | L | Both | R |
|-----|---|------|---|
|     | (17+18(1))(5) | (17+18(1))(5+18(1)) | (17)(5+18(1)) |
|     | (35)(5)=175 | (35)(23)=805 | (17)(23)=391 |
|     | (LOPI 17)(LOPI 5)=LOPI 13 | (LOPI 17)(LOPI 5)=LOPI 13 | (LOPI 17)(LOPI 5)=LOPI 13 |
|     | 175=13(mod 18) | 805=13(mod 18) | 391=13(mod 18) |

Essentially each seed operation gives three options each time, L=1/R=0, L=1/R=1, L=0/R=1, as shown above and in the ternary diagram. Each new number generated is treated in the same way. Given that some of the values begin to repeat that are not due to unique factors that would interfere with the calculation of the number of factors per integer, as discussed later, the best way is linear binomial convolution of the two sequences with the equation, $(LOPI_1 + 18(a)(LOPI_2 + 18(b)))= N$ LOPI

Therefore, the binomial convolution for subset generated from the seed 17x5 is

$((17 + 18(a))((5+ 18(b))$, where $a \geq 0, b \geq 0$.

This can also be written in the polynomial form;

$85 + 306(b) + 90(a) + 324(ab)$,

where f(a,b)=N LOPI element, only for non-primes and therefore the gaps in the range of the codomain are the primes left in each Superset after convolution. This gives us a good piece of information to use to prove that the primes are also produced at a cyclic rate as well, as the gaps are cyclic since the Non-prime codomain is cyclic. In the polynomial form however, generation of elements of the subsets is less efficient and the factors are not readily visible, only the solution, although random coefficients can be put in to generate different random non-prime numbers. Subset generation reveals the factors being multiplied to produce the elements within the LOPI NP subsets. An interesting aspect about subset population is the idea that we have a set of subsets that are always in the NP state and each of these growa as we increase a, and/or b of the convolution.



It would be interesting to model the behavior of the changing state of the total LOPI Superset as a whole.  Sums and total number of primes and non-primes as they change in time, as NP creates a new element, it is subtracted from the Superset, but not until a later cycle in the linear Superset.  The state of the system is constantly changing and here we could understand an equilibrium, a distribution and a consequence to this "equilibrium".   A constantly growing Superset, that is preceded by a constantly growing NP subset, that creates a new space in which future elements in the linear equation get generated first into the planar and solid numbers. The average at any point of the linear progression appears to be equal to the value of the LOPI for the reduced residue class mod 18.

Each LOPI has also a subset generated from itself with the requirement that $b>0$, since each LOPI is prime except LOPI 1. For LOPI 13 and seed operation $13 \times 1 = 13$, we must, by definition, include this as a possible subset since this seed multiplication does yield 13, although it is the LOPI identity for 13.

$$13 \times 1 \text{ becomes}$$

$$(13 + 18(a))(1 + 18(b)), a\geq 0, b>0$$

thereby forbidding multiplication by the trivial factor of 1 and blocking the mapping of primes into the NP subsets by linear convolution.  This leaves all primes partitioned out and left behind as the unique elements of the original LOPI superset.

**See Figure 3**

$$(13+18)(1) \quad (13+18)(1+18) \quad (13)(1+18)$$

The Left branch equals the superset, which contains all elements of LOPI 13, Prime and Non-prime.  We avoid this in the convolution by not allowing b=0, as already discussed, thereby also avoiding the prime 13 from being in the NP subset.  The linear convolution for r.r. seeds LOPI 13,  13 and 1: $((13+18(a))((1+18(b)), a\geq 0, b>0$. Another note about subsets made from r.r. totatives of 18 that are equal, or 7x7, to avoid the unnecessary multiplicity the coefficients a and b are $a\geq b\leq a$.

Some LOPI classes have a square as one of its seed operations such as LOPI 13 seed 7x7. When this occurs, the polynomial from the convolution of the two binomial linear equations changes to define a in terms of b, with $a \leq b \geq 0$. From the number line, figure 1, in this manner building from the first unique lowest non-primes of the seed operations, the binomial linear equations for convolution to produce the data matrix are evident. All the seed operations of each of LOPI are listed below along with the generating binomials for the data function table.

The r.r. combinations for the NP generating linear binomial convolutions in each LOPI set yield:

**LOPI 1**:                                              LOPI    NP    Polynomial: LOPI element + Neutrals=$N_x$ LOPI

1x1 =    $((1 + 18(a))x((1+18(b))$,    $1 \leq a \leq b \geq 1$      $1 + 18(a+b) + 324(ab)$

7x13=    $((7 + 18(a))x((13 + 18(b))$, $a \geq 0, b \geq 0$      $91 + 126(b) + 234(a) + 324(ab)$

11x5    $((11 + 18(a))x((5 + 18(b))$, $a \geq 0, b \geq 0$      $55 + 198(b) + 90(a) + 324(ab)$

17x17    $((17 + 18(a))x((17 + 18(b))$, $0 \leq a \leq b \geq 0$      $289 + 306(a+b) + 324(ab)$

**LOPI 5**:

5x1    $((5 + 18(a))x((1 + 18(b))$, $a \geq 0, b > 0$      $5 + 90(b) + 18(a) + 324(ab)$

11x7    $((11 + 18(a))x((7 + 18(b))$, $a \geq 0, b \geq 0$      $77 + 198(b) + 126(a) + 324(ab)$

13x17    $((13 + 18(a))x((17 + 18(b))$, $a \geq 0, b \geq 0$      $221 + 234(b) + 306(a) + 324(ab)$

**LOPI 7**:

7x1    $((7 + 18(a))x((1 + 18(b))$, $a \geq 0, b > 0$      $7 + 126(b) + 18(a) + 324(ab)$

5x5    $((5 + 18(a))x((5 + 18(b))$, $0 \leq a \leq b \geq 0$      $25 + 90(a+b) + 324(ab)$

13x13    $((13 + 18(a))x((13 + 18(b))$, $0 \leq a \leq b \geq 0$      $169 + 234(a+b) + 324(ab)$

11x17    $((11 + 18(a))x((17 + 18(b))$, $a \geq 0, b \geq 0$      $187 + 198(b) + 306(a) + 324(ab)$

**LOPI 11**:

| | | |
|---|---|---|
| 11x1 | $((11 + 18(a))\times((1 + 18(b))$, $a\geq 0, b>0$ | $11 + 198(b) + 18(a) + 324(ab)$ |
| 5x13 | $((5 + 18(a))\times((13 + 18(b))$, $a\geq 0, b\geq 0$ | $65 + 90(b) + 234(a) + 324(ab)$ |
| 7x17 | $((7 + 18(a))\times((17 + 18(b))$, $a\geq 0, b\geq 0$ | $119 + 126(b) + 306(a) + 324(ab)$ |

**LOPI 13**:

| | | |
|---|---|---|
| 13x1 | $((13 + 18(a))\times((1 + 18(b))$, $a\geq 0, b>0$ | $13 + 234(b) + 18(a) + 324(ab)$ |
| 7x7 | $((7 + 18(a))\times((7 + 18(b))$, $0\leq a\leq b\geq 0$ | $49 + 126(a+b) + 324(ab)$ |
| 11x11 | $((11 + 18(a))\times((11 + 18(b))$, $0\leq a\leq b\geq 0$ | $121 + 198(a+b) + 324(ab)$ |
| 17x5 | $((17 + 18(a))\times((5 + 18(b))$, $a\geq 0, b\geq 0$ | $85 + 306(b) + 90(a) + 324(ab)$ |

**LOPI 17:**   Contains all twin secondary primes when a=b,

| | | |
|---|---|---|
| 17x1 | $((17 + 18(a))\times((1 + 18(b))$, $a\geq 0, b>0$ | $17 + 306(b) + 18(a) + 324(ab)$ |
| 5x7 | $((5 + 18(a))\times((7 + 18(b))$, $a\geq 0, b\geq 0$ | $35 + 90(b) + 126(a) + 324(ab)$ |
| 11x13 | $((11 + 18(a))\times((13 + 18(b))$, $a\geq 0, b\geq 0$ | $143 + 198(b) + 234(a) + 324(ab)$ |

{U}=S1{1+18($a_0$), (1 + 18($a_0$+1), (1+18($a_1$ +1)...}, the same for all S5 {5+18(a)...}, S7 {7+18(a)...}, S11{11+18(a)...}, S13{13+18(a)...} S17{17+18(a)...}, where a≥0.

We have work backwards in a sense, in that we group the integers by an equivalence identity, the core r.r. LOPI, and found that by recognizing that all r.r. factors in the lowest NP seed products are the LOPI reduced residues of mod 18 that will go into the linear convolution for each LOPI congruence class and will continue to produce all the non-primes cyclically through convolution, while maintaining the core LOPI grouping identity within the NP subsets to infinity. The two linear equations in the convolution are the rings of divisors for each subset within each LOPI superset all coprime to mod 18. The non-primes in each LOPI class are reducible by the non-trivial factors generated in each polynomial pair in the Dirichlet linear convolution,

(NP LOPI + 18(c))/(LOPI$_1$ + 18(a) = LOPI$_2$ + 18(b)

and are each LOPI cores that return upon multiplication the product with the same LOPI core of the original grouping identity, the LOPI residue mod 18.

Each element in the LOPI congruence classes, f(a) in the linear Dirichlet polynomial, primes are the only elements with trivial factorization and therefore are also the elements for which the polynomial is not irreducible f(a)=LOPI + 18(a) in the natural numbers ((g.c.d. N,18)=1). For prime elements in N no two linear progressions reduced residue mod 18 exist that can multiply to yield that result, unless we consider holding the second linear equation in ((LOPI + 18(a))( LOPI 1 + 18(b)), where b=0, therefore we get right back to the Dirichlet linear progression that created the LOPI congruence classes originally, LOPI + 18(a)(1 + 18(b)) with b at zero, equals LOPI + 18(a)(1) or the original Dirichlet LOPI mod 18 progressions.

The elements within the subsets can be ranked in terms of primality based on their multiplicity within the subsets for LOPI class in which they are found by (2(number of times found) + 2(trivial))/2, with Prime 2(0) + 2=2/2=1, they will not be found in the subsets. We note here one special case on multiplicity of the squares when multiplied by themselves, as discussed on page 4. In the LOPI$^2$(LOPI) the product will look like a secondary prime due to the fact that its factors are 7X49, or 7X7X7, in these cases only where the prime factorization is the LOPI$^3$ show a multiplicity of one in their respective LOPI subset. This is the only multiplicity of one that looks like a secondary prime but is a square secondary prime and therefore of low order and can be easily understood. Return to page 4 for the full discussion of this unique false positive of multiplicity 1 for a secondary prime, noted again here as: 5^3, 125, 7$^3$ 343, 11$^3$ 1331, 13$^3$ 2197,17$^3$ 4913.

We can define the primality of an integer in two ways given the data function factor/product tables from the convolution of the unique binomials. If we are concerned only in the number of divisors of an odd natural number, then we define (d)divisors=#of times found in subset of LOPI superset * 2 + 2 for the trivial, or we can ask if a number is prime by querying the number of times a natural odd integer occurs in a subset, if 0 we know we have a prime, with only the 2 trivial divisors, 2*0+2=2, 0 for number of times found in its particular LOPI subsets. The secondary prime is found in its unique LOPI superset and in 1 non-prime subset, 1 time, Secondary prime = 2(1) +2=4 total factors including the two trivial factors, self and one. Generally, we can write:

Element N multiplicity = (2a)=degree of primality, where a = number of times upon query that element N occurs in the specific proper non-prime subset for the LOPI N.

In a Venn diagram of figure 2, the intersection of the subsets, a prime is not found in the subset and a secondary prime is not found in the intersection of any of the subsets. Secondary primes are uniquely within only the seed subset, no intersection to any other seed subset within the LOPI N occurs.

In this article we present the basis of a new system of grouping all natural odd numbers through a reduced residue mod 18 constant partition identity and show how the factors of all the elements in the non-prime subsets are cyclically and deterministically generated through linear convolution, while the primes within each congruence class are generated in time as the gaps in the convolution of the two binomial linear equations. Both binomials are in the form of Dirichlet progressions, LOPI +18(a), with LOPI and 18 coprime, therefore upon convolution each subset contains an infinite number of non-primes. Each LOPI congruence class used as a factor generator in the Dirichlet linear binomial convolution, $(LOPI_1 + 18(a))((LOPI_2 + 18(b))$, contains an infinitude of primes as well. Once the subsets are subtracted from the original superset we are left with a partitioned set of primes for each r.r. LOPI mod 18 family, 6 partitioned subsets of primes, that compose the whole of the Prime set within the Universal set of odd integers coprime to mod 18. According to Dirichlet the number of primes is infinite in a linear progression, LOPI X + B(q), where g.c.d. LOPI X,B= 1, therefore we have partitioned out the prime elements of the 6 infinite LOPI congruence classes foreach LOPI N class and when added together:

{P LOPI 1 }+ {P LOPI 5} +{LOPI7}+{LOPI 11}+{LOPI 13}+{LOPI 17} = $\mathbb{P}$

The total set of Primes in our Universal set, the Natural number set of all odd integers coprime to 18, which covers all possible prime integers, excluding the outliers of 2 and 3 as discussed earlier.

The subsets generate non-primes more quickly since there are 3-4 linear binomials being convoluted while only one binomial arithmetic sequence of LOPI + 18(a), populates the congruence class of the LOPI N superset. In the 3-4 generating subset functions which generate all the non-primes for each LOPISuperset, the Prime set left after the subsets are subtracted out is clean and complete. The LOPI mod 18 congruence class systems leave all primes gathered in 6 distinct sets organized by their unique LOPI identity which equals the residue mod 18, $N_{LOPI} \equiv LOPI(mod\,18)$, which were orignally generated into the LOPI superset by the linear progression, LOPI + 18(a), with a≥0.

The primes are present and there are no missing non-prime elements upon generation of the factors by point wise convolution into the subsets, therefore integer multiplicity in the subsets is a valid deterministic method by LOPI mod 18 for primality testing and factor generation. A short example, if we look at the non-primes generated by the binomial expression of LOPI 13, we see that the non-primes in the data function matrix of each non-prime that exists within each set and the factors are in the rows and columns of each input, row -1 = a, and column - 1, equals b.

LOPI 13, the all superset elements up to 301:

13, 31, 49, 67, 85, 103, 121, 139, 157, 175, 193, 211, 229, 247, 265, 283 e 301.

List the non-primes generated by the seed binomial LOPIs:

13 +18(a)(1+ 18(b)= 247, for our analysis up to 301 as mentioned above, we go no further since there is no intersection here in this limit.

7x7;( (7 +18(a))((7 +18(b)); a≤b≥0 : 49, 175, 301

11x11; (11+18(a))((11+18(b) where a≤b≥0; : 121;        7 total non-prime from the subsets

17x5;(17+18(a)(5+18(b); a,b≥0 : 85, 175, 265.

With the data matrix subset we query and find 301 is generated from the seed operation for LOPI 13 by 7x7, or (7 +18(0)(7 + 18(2)).

 301 is in row 1, column 3, or 1-1=0=a, factor 7 and column 3, 3-1=2=b or (7+ 18(0)(7+18(2), which gives the factors 43 and 7, 43(7)=301.

The primes left are 13, 31, 67, 103, 139, 157, 193, 211, 229, and 283.

As a prime counting function we see that 301 is $X_a$=16, $X_n$=16+1=17, we have 17 elements including 301. We generated 10 primes that weren't in the subsets therefore we should be able to account for the 7 non-primes from the subsets, which are listed above and do equal 7 in total. All 17 elements including 301, are generated and accounted for. By the Prime counting function, N/lnN, the estimate is about 53 total, but we do not count 2,3, therefore 53-2=51/6 to put into the 6 LOPI classes, puts us at about 9 per class, which is reflecting an equal distribution overall into the 6 LOPI congruence classes mod 18.

Another example is included later with a delta function to show the growth in the subsets to show again the distribution averages out between residue classes as Dirichlet and the PNT suggest. If we didn't have each subset generating the factors and the product of the other non-prime odd natural numbers, we would have mistaken some numbers as prime. The binomial linear Dirichlet progressions used in mod 18 with the LOPI as constants, are like gears running in different times and when coordinated through the subset convolution, the primes in each LOPI are the numbers **not** generated. They cannot be generated as the factors that are generated by the "gears" of the subsets never multiply to yield primes since their only LOPI r.r. factor mod 18 of 1, is never used, and therefore, they are skipped in the cycle at which the "gears" change. When we drive we often skip a gear while going from $3^{rd}$ to $5^{th}$ when the speed is right, skipping the prime's gear is akin to that, we skip convolution by LOPI 1, thereby leaving the primes in the LOPI superset, being only generated linearly by their linear Dirichlet expression LOPI + 18(a), a≥0.

Any number, therefore, can be queried for its multiplicity within the proper non-prime subset matrices and the primes yield: found=0

and the unique secondary primes yield a find of 1.

The unique factors of each queried integer are known as well. Figure 4 shows more of the non-prime convolution data matrix for LOPI 17.

For example 6313, query in the specific range listed in the attached algorithm yields:

 output: found 1 time  - subset 17*5, row 6, column 4

We know then that this is a secondary prime, found one time, 2(1) + 2(trivial)=4, total, 2 trivial and non-trivial factors are (17 +18(row-1)(5+18(column-1);

 (17 +18(5)(5+18(3)=6313= 107(59)

If the four queries into the subsets had been output: not found, then n=0, and our integer would have been a prime, left partitioned out only in LOPI N Superset and with 2 divisors, self and 1 by 2+2(0) =2, a prime.

0 finds in specific LOPI subsets query = prime.

Each LOPI superset and its subsets go to infinity. The Dirichlet arithmetic progression contains an infinitude of primes by the Linear Binomials, LOPI + 18(a ), $a≥0$, and a={N}, which also extends to infinity, we therefore know that the non-primes in each LOPI continue to infinity as well, prime and non-prime distribution is even throughout each LOPI congruence class.

With the LOPI system of binomial convolutions, we are generating the non-primes and their factors, rather than inquiring as to whether a number is divisible by certain numbers or how many factors. Primality = number of subsets in which we find the integer, 0, being prime, 1 is a secondary prime, 2 a tertiary and so forth as discussed above with the relationship, 2 + 2(n), with n=numbers of times number(n) is found in a subset, the result equals the number of factors of n, with primes never being found in a linear convolution subset, therefore n=0 and primality=0, and number of factors = 2, trivial. For the negative integers, we cross over left to zero and the lowest odd partition integer becomes equivalent to a positive counterpart -18. For example, 5-18=-13. The LOPI 5 + 18(a) becomes the (– 13 residue class upon crossing the y axis in the negative direction, (-13 LOPI + 18(-a), with 1 becoming -17, 7 to -11, these are found to be the inverses for the LOPI mod 18, |G| 6, (1,5,7,11,13,17) finite cyclic group of residue elements generated by $5^g$, $11^g$ with a repeat every 6, as $5^{1+6(a)}$, $5^{3+6(a)}, 5^{4+6(a)} 5^{5+6(a)} 5^{6+6(a)},$, also with LOPI 11. **[12, 17]**. Generators 5, 11 are found to generate the LOPI residue pattern, not all the actual integers, but the LOPI residue classes we use in the Dirichlet linear progressions that lead to the elements of all the congruence classes of the odd prime/non-prime integers upon iteration of 18. For example the first iteration of LOPI 7 is generated by (7+18(1)= $5^2$=25=$X_{a=1}$ and we can cyclically see LOPI 7 elements in sets of $5^{2+6(a)}$, with a≥0, therefore an exponentiation cyclic LOPI 7 element again at $5^8$=390625 =LOPI 7 Element number X $_{a=21701}$. $5^8$-$5^2$=390600 which indicates that at cycle 390600 of 18 or 21701-1 * 18 = 390600 +(18 +7=25)=390625 we have element 21702, $X_n$, as every iteration of 18 yields a new cycle with a new

LOPI element upon (LOPI + 18(a)).  We can know a larger range of LOPI class elements in a larger cyclic non-prime pattern from within which to do different pattern analysis. The generator 11 perhaps is also useful in another way, as any LOPI element ending in 5 is easy to spot as non-prime.  LOPI generator 11 gives us less obvious non-primes created cyclically in the LOPI 7 congruence class by  $11^{4+6(a), a\geq 0}$ that are all multiplies of 11, therefore we know that all large elements in LOPI 7 congruence class generated by $11^{4+6(a)}$ are non-prime and divisible by 11 and by all the elements $11^{\leq 4+6(a)}$ . 11 as a residue cyclic generator mod 18 yields very useful information about these larger elements in each LOPI residue class as to their primality and their factors without searching the data matrix.  We can say that all elements LOPI 7 produced by these cyclic residue generators mod 18 are non-prime and they are divisible by the generator $11^{4+6(80)}$.  LOPI 7 elements generated by $11^{4+6(a)}$ are always non-prime, therefore we can take:

 log(LOPI 7 integer)/log(11)

to see that the large element in a specific LOPI  7 is a non-prime divisible by 11 or log(LOPI element)/LOPI 11 yields exponent in $11^x$, thereby knowing the LOPI class as well, as listed above.     This can be a deterministic test for primality and factorization for some larger integers that may fall in the cycle of residue generation via LOPI $11^{1-5\ +6(a)}$.  Ultimately the convolution/query of the specific data matrix tables in the integral cycle needed is more efficient and creates all possibilities, but to make a list is possible as well for at least the elements that fall within these parameters by LOPI 5, LOPI 11 as cyclic residue element generators.  This gives us an analytical tool to apply in different situations, for instance in looking for Mersenne primes.  The LOPI 11 elements fall right after all the perfect numbers, which are all LEPI 10, as will be discussed later.  If we are able to identify the perfect number as a Euler totient, then we know that the following LOPI 11 element is a prime and therefore should not be in the data matrix table and should not be log(element)/log 11 = exp.1+6(a).

The set of LOPI congruence classes with their elements are:

{U r.r. mod 18}=S1{1+18(a)...}, S5 {5+18(a)...}, S7 {7+18(a)…}, S11{11+18(a)...}, S13{13+18(a)...} S17{17+18(a)…}, where a≥0.  All the subsets generated by the seed operations in each LOPI are proper subsets of A**,** or A is a proper superset of B. The non-primes all generate from the convolution of one of the r.r. factor seed operations used in the linear binomials, for example 17(5)=85=LOPI 13, therefore the convolution ((17 + 18(a))((5 + 18(b)), a≥0,b≥0, will generate all the non-prime elements and their factors for this NP subset of LOPI13, Subset B, Figure 2.

All the seed operations are treated in this manner to generate all the elements of the other subsets that exist within each LOPI. Superset LOPI progression produced by the LOPI +18(a), with a≥0, not to be confused with the non-prime subset linear convolution ((LOPI $_1$+18(a)(1 + 18(b)), b>0.

 (1+18(b)) needs to be discussed here as it's different, where LOPI 1 and 2 are the same plus they are LOPI$^{\wedge 2}$, therefore to avoid repetition because LOPI 1 and 2 are equal we put that a≤b, and to avoid any multiplication by 1 we must say 1≤a≤b≥1, while the Superset linear Dirichlet binomial is LOPI 1 is 1 + 18(a) a≥0, as all the other LOPI Supersets. All subsets of LOPI 5,7,11,13,17 +18(a)(1+18(b)), with b>0, a≥0, where 1 is LOPI$_2$ b must never be zero, thereby avoiding that prime elements in the LOPI supersets be partitioned into a non-prime subset. When LOPI 1 and LOPI 2 are equal, other than LOPI 1 as already discussed, or, LOPI$_1$+18(a)(LOPI$_2$+18(b), where LOPI$_1$=LOPI$_2$ we use a≤b, with b≥0. For example for the partitioning of LOPI 13 into its non-prime elements from the e.r. 7*7=49=LOPI 13, (7+18(a)(7+18(b))=LOPI 13 N nonprime subset elements when we set a≤b and **b≥0.** These criteria are few and serve to keep 1 out as a factor which serves to keep the primes in the LOPI superset as they have only r.r. mod 18 of 1, and the non-primes partitioned out into their respective LOPI non-prime subsets. All criteria of a and b are listed above with each subset for the 6 LOPI congruence classes. Prime numbers cannot be generated from the convolution of these binomials and are formed only from the original Dirichlet linear progression of LOPI + 18(a), with a≥0. The pattern in the primes depends upon the frequency of the occurrence of a new factor from the seed operations as discussed above. From our point of view, it is essential to see the numbers in their LOPI groupings and the rhythm of the factors of the non-primes within each subset of each LOPI superset. The primes can be partitioned out and highlighted due to the pattern of the non-primes and their factor generation mechanism being based on their Lowest Odd Partition Integer, the primes are generated as well, but only by multiplication by LOPI 1 and we limit our system in this way, no multiplication by one is allowed into the non-prime LOPI subsets. Primes can not be mapped into the NP subsets and are left behind to remain within their LOPI superset, as all non-primes are cyclically generated and partitioned out into their unique LOPI NP subsets.

On the ternary tree, figure 3, the non-primes can be seen to generated at large intervals as they are due to the convolution of two binomials, whereas the primes and non-primes are generated step wise by +18 in the linear progression, therefore the primes are easily recognized in time within each LOPI group and can be partitioned out, as their location is well defined in contrast by the convolution of the values of non-primes. The LOPI congruence classes containing all possible odd natural number primes are generated by one Euclidean, Dirichlet equation, the non-primes are generated by the gears of their seed convolutions, skipping the primes, allowing for prime identification due to a lack of presence in the NP subset elements generated by the convolution and therefore subsequent multiplicity of 0 in all proper non-prime subsets of the respective LOPI of the integer. Primes are readily identifiable, as the equivalence relation does not map the LOPI primes in A to the LOPI proper subsets of non-primes, there are no primes in any of the LOPI subsets, as multiplication by 1 is not allowed and there are no non primes left in the LOPI superset after convolution, as well the linear convolution that generates all the non-prime subset elements well outruns the linear progression populating the LOPI Super-sets, allowing steps in time to partition out the primes as well at whatever upper value we decide. We have a partitioned subset of primes in each LOPI congruence class to work with as we wish, for instance we can multiply them to make new and known large secondary primes, or identify in efficient time if an integer is a Sophie or Mersenne prime, or a member of a twin prime set.

In conjunction with the data function matrices generated for all the elements of the non-primes in each LOPI congruence class mod 18 we present the attached arithmetic algorithm as a deterministic efficient primality test based on unique multiplicity, where a prime returns a 0 upon query, as primes are not present in any of the non-prime subsets of each LOPI. As well, the algorithm reveals deterministically and in efficient big(0) time, (see Appendix), all the non-trivial factors of all the non-primes in each LOPI subset by convolution of the linear binomial progressions derived from the lowest non-prime seed factors in each LOPI linear progression mod 18, as shown in Figure 1. In the next sections of this work we analyze various types of primes using the LOPI reduced residue system mod 18.

**Sophie Germain Primes**

The 6 LOPI super-sets and the rules for the generation of their elements have implications to the possibility of Sophie Germain Primes **[18]** co-occurring as twin primes and when a Sophie Germain Prime can exist within a LOPI superset. Our work has found that no Sophie Germain Prime can exist in LOPIs 1, 7, 13, $N_{LOPI} \equiv 1 (mod 18)$, $N \equiv 7 (mod 18)$, $N \equiv 13 Mod(18)$, which can be generated together in one average inclusive equivalence relation for all elements in our universal set of congruence classes coprime, odd natural numbers mod 18, $N \equiv 1 (mod 6)$. No Sophie Germain Primes $\equiv 1 (mod 6)$ holds. LOPI congruence classes 1, 7, 13 will have no Sophie Primes in their Dirichlet progressions, ($LOPI_{1,7,13} + 18(a)$), or odd natural numbers N coprime to mod 18, $N \equiv 1 (mod 6)$, cannot be Sophie Germain primes. On the other hand, the other twin LOPIs to these three, LOPI 5, 11, 17, have the possibility of being a Sophie prime, all are generated under the average congruence relation of $N \equiv 5 (mod 6)$. Given that only LOPI 5, 11, 17 elements can only be Sophie primes, the equivalence relation check of $3 (mod 4)$ derived from Gauss' Law of Quadratic Reciprocity **[19]** used to find a matching safe prime of 2p + 1 to be a divisor for a Mersenne number $2^p - 1$, can be limited only to the LOPI congruence classes of 5,11,17. This is a more exact reduction since in the non-possible LOPI classes, 1,7 and 13, we see also the equivalence relations $3 (mod 4)$ met for elements that cannot possibly be Sophie Germain primes by the relationships discussed above, but become justifiably exempt under the congruence relation $5 (mod 6)$. For example, 79 is congruent to $3 (mod 4)$, but it cannot possibly be a Sophie Germain Prime as it is in LOPI class 7 and upon the operation 2p +1 will, as shown above, always yield a LOPI 15 element, which we know can never be prime and are not generated in any LOPI residue class capable of containing odd natural primes. We work with all residue classes coprime to 18. Based on the reduced residues mod 18, that also equal the core LOPI that also equal the totatives of 18 and residue 15 is not a reduced residue mod 18. Computation time can be cut dramatically by automatically taking out from the search any of the elements in LOPI classes 1,7,13 as by current methods these elements may meet the 3(mod 4) relation and may be checked to see if they meet the Sophie prime requirement before testing them as divisors of Mersenne numbers, but this isn't necessary, any element meeting the equivalence relation 1(mod 6) should not be pursued as a potential Sophie Germain prime. Only elements meeting the relation 5(mod 6) will yield a possible Sophie Germain safe prime, elements in LOPI 5,11,17 mod 18 congruence classes. One could also just do the $N_{LOPI} \equiv LOPI(mod\ 18)$, or add the sum digit choosing the odd equivalent of the partition identity and also therefore know not to purse any integer that has a LOPI of 1, 7 and 13.

Modular math follows the same rules as regular math in that an integer element from the LOPI congruence class undergoing the Sophie Germain operation of $((LOPI +18(a)) * 2) +1=N$ LOPI is predictable as to its internal core lowest odd partition identity and therefore its LOPI congruence class and cannot change, for instance;

$((7 +18(a))*2)+1=15 +18(a)$,    $((1+18(a))*2)+1=3 + 18(a)$, and    $((13 + 18(a))*2)+1=9 +18(a)$.

These residue elements mod 18 will always generate elements in the congruence classes $3 + 18(a)$, $15 + 18(a)$ and $9+18(a)$, with $a \geq 0$ when we carry out the Sophie Germain operation onto the elements in LOPI classes 1, 7, 13. Therefore LOPI classes 3,15 and 9 as discussed above, never contain primes as they are not coprime totatives to $\varphi(18)$, equivalently they are not LOPI identities that can yield primes, under $LOPI + 18(a)$, with $a \geq 0$. No Sophie Germain Safe Primes are found in LOPI r.r. 1, 7 and 13 mod 18 residue classes as upon the operation we just generate elements for the not coprime residues 3, 9 and 15 to mod 18. Equivalently we can say that any element from the LOPI congruence classes 1,7 and 13 will never reveal a prime by the inverse of the Sophie operation, $(N\ LOPI_{1,7,13}-1)/2$.

Upon analysis using the details of LOPI mod 18, we can see that no twin prime pair of primes can simultaneously both be Sophie primes as well, one or the other may be a Sophie Germain prime, but we cannot have two Sophie Germain primes who are two twin primes, except with the one time exception of integer 3, since 3 is the unique and only prime within the congruence class 3, and 3(2) +1=7, a prime and the twin of 3 is 5, which is a Sophie Germain prime as well. We therefore have the unique exception twin set 3,5 as two twin primes that are both Sophie Germain primes. This is an exception based on the outlier 3 as a prime, but not a reduced residue of 18. The integer 3, however does not meet the inverse action of the Sophie Germain safe prime operation, as 3-1/2=1 and 1 is not considered prime and in general a Sophie Germain safe prime is always preceded by a prime under the operation, such as 107-1/2=53, therefore 53 is a Sophie Germain with the safe prime of 107. Integer 3, however, is considered a Sophie Germain prime since 3*2 +1=7, but 3-1/2=1, where 1 is not prime, so in reality to meet the definition it appears that we either have to reconsider 3 as a non-prime or reconsider 1 as a prime.

We can see that all twin prime LOPI pairs, [1,17], [5,7],[11,13] have one of each type of LOPI class elements N≡1(mod 6) for (LOPI 1,7,13(mod 18), and the twin is of the LOPI classes, N≡5(mod 6) for (5,11,17(mod 18). This grouping information is helpful if searching for large Sophie Germain Primes that may also be twins, as in the research in 1999, of Indlekofer and Jarai [20], which indicates that they were also searching for a large Sophie Germain prime pair. Based on the analysis above we see now that this is not possible and why their search did not yield this result. When we have a pair of twins that are prime, they will never have the simultaneous relationship of being Sophie Germain Primes and also Twin Primes. Here again we see support for the even distribution of the infinite primes contained in the reduced residue system of ɷ(6), [1,5], **N≡1(mod 6), N≡5(mod 6),** as these two congruence relationships separate all possible odd primes into two larger disjoint congruence classes, each containing the one of the twins in a twin prime set, by the definition of twin prime= P, P+2. Mod 18 is a finer partitioning of this broader average equivalence relation used here to discuss Sophie Germain primes and can be used to save an important amount of computational energy and time by limiting the search to LOPI r.r. 5, 11 and 17 congruence classes mod 18.

Again we refer to LOPI 7, 7*2 + 1= 15, LOPI 15 will be the LOPI congruence class for LOPI group 7 elements that undergo the multiplication by 2 and addition of 1, as the operation is essentially on the core LOPI value, as in this LOPI mod 18 system each integer is viewed as an unrestricted partition of the lowest odd partition identity, or LOPI mod 18. An integer no matter how large can be separated as its core and a quantity of cycles of 18, a neutral quantity added to the core LOPI, thereby maintaining the core LOPI identity through addition of a neutral 18(a), LOPI + 18(a).   Sophie Germain Prime integers exist in LOPI congruence classes 5, 11 and 17 and no twin prime pair exists that are both Sophie primes, since twin primes are LOPI 5,7;  LOPI 11,13; and LOPI 1,17, and each set has only one LOPI that can be a Sophie Germain prime. The LOPI congruence classes can be represented as a larger congruence relation, LOPI 1, 7, 13 $\equiv 1(\bmod 6)$ and LOPI 5, 11, 17$\equiv 5(\bmod 6)$, the latter representing the classes in which we may find a Sophie Germain prime.  The largest Sophie Germain prime will not have an accompanying twin prime, as deterministically proven by LOPI mod 18. Sophie Germain Primes are only in the congruence relation N$\equiv 5(\bmod 6)$, or LOPI 5, 11, 17.

 Interestingly also the LOPI 5, 11, 17 are as well residue cores that are never formed by squaring a residue, only LOPI class elements 1,7,13 are generated upon using 2 as the exponent as discussed in the generator for the residue pattern. Therefore we can extend our research in the future to the question of Fermat's Last Theorem, that no exponent other than 2 can form a perfect square under Pythagorean's Theorem, $a^2 + b^2 = c^2$.  We know that a perfect square can not be found in LOPI congruence classes 5, 11 and 17, which gives us a combinatorial restriction for $a^2 + b^2$, that it can not add to yield an element within the Sophie Germain possible LOPI classes of 5, 11 and 17.  It appears that we must always have a neutral element on one side of the equation or the other to produce a Pythagorean Triple, or $c^2$ is a perfect square from disjoint LOPI residue classes 1, 7 or 13, as well we see the perfect square can also be from the neutral non LOPI residue classes 3, 9 and 15 mod 18. This is a future analysis to be completed.

**Mersenne Primes**

This discussion leads naturally into the discussion of possible Mersenne numbers and Mersenne primes **[21]** as they are found to be in the other set of LOPI congruence classes, LOPI 1, 7, 13$\equiv 1(\bmod 6)$.

Interestingly there is only one Mersenne prime that can be in the LOPI residue class 7, encountering here the outlier in the LOPI 3 congruence class.  As mentioned earlier the totatives of 18 are the reduced residue system for mod 18 and they match exactly the value of the possible lowest odd prime identity integer, but 3 is not a totative of 18 and cannot produce a superset of elements mixed with odd prime and non-prime.  It produces only non-prime elements, except for the LOPI itself, 3, which is considered prime. Therefore we make a note here that $2^3-1=7$, and that 3 as a LOPI to generate more exponents for possible Mersenne prime generation does not work as they are all divisible by 7, the original unique prime Mersenne generated by the unique LOPI 3 prime element, which is not a coprime residue mod 18, $2^{3+6(a)} - 1$ all produce LOPI congruence class 7 elements divisible by 7, or LOPI 7 non-prime elements that meet the congruence relation N≡1(mod 6). LOPI 7 elements cannot contain Mersenne prime elements, nor can they form Sophie Germain primes, but the can form Mersenne primes by $2^{1+6(a),\ or\ Non\text{-}S.G.\ class\ primes}$. For example $2^9 - 1=511$, element LOPI 7, 511 is congruent to 1(mod 6), therefore we know it is not a Sophie Germain and on the flip side it isn't congruent to 5(mod 6), therefore not a Sophie Germain prime. However if we only used the 3(mod 4) equivalence relation, we may say that yes this can divide a Mersenne number, it is perhaps a Sophie Prime, but if we use 5(mod 6) as a requirement for Sophie Germain inclusion then we would know immediately that 511 is not a Sophie Germain Prime, also from the beginning we can see that 511 is a LOPI 7 element and can not be a Sophie Germain prime, safe or unsafe, nor can it be a Mersenne prime.

The only Mersenne primes we can generate are going to be found in LOPI 1, and LOPI 13; LOPI 1 Mersenne being generated by the exponent elements of $2^{LOPI\ 1\ +\ 6(a)}$, LOPI 1,7,13 exponent classes and $2^{LOPI\ 5\ +\ 6(a)}$ generating the LOPI 13 Mersenne primes, p=primes from the LOPI 5,11,17, (the Sophie Germain Primes). Mersenne primes located in LOPI 1 will have been generated by exponents in the LOPI congruence classes of 1,7,13, while Mersenne primes in LOPI 13 will have been generated by exponents in the LOPI congruence classes possible to contain Sophie Germain Primes, LOPI 5,11,17, the twin LOPI congruence classes of LOPI 1,7,13. We can therefore conclude that any Mersenne prime will not be a Sophie Germain prime as they exist in disjoint reduced residue classes mod 18. Mersenne Primes have a LOPI core of 1, 13 as generated by exponentiation of either the Sophie Germain primes, $2^{LOPI\ 5\ +\ 6(a}$, p=primes from the LOPI 5,11,17 classes, N≡5(mod 6), yielding Mersenne prime elements in LOPI 13 class, and Mersenne prime elements in LOPI class 1, generated by the non-Sophie Germain classes, $2^{LOPI\ 1+\ 6(a)}$, representing the primes from the LOPI 1,7,13 classes, N≡1(mod 6). By partitioning the Mersennes primes into Sophie Prime or Non-Sophie prime generating primes we can know the class the Mersenne prime generated from depending on its LOPI congruence class. For example Mersenne 131071, without doing any calculation except using the LOPI/LEPI core identity we know that Mersenne 131071 is a LOPI 13 element and therefore originated from a Sophie Germain prime. We know that Log 131071+1/Log 2 = a prime in the LOPI classes LOPI 5, 13 or 17. Working backwards we see that log 131071/log 2= 17 a possible Sophie Germain prime exponent in LOPI 17 class. Also for the next Mersenne prime 524287, LOPI 1, we know that this exponent is an element of the Non-Sophie primes classes, 1, 7, 13 and in fact find that 52487+1=log 52488/log 2=19, a Non-Sophie Germain prime.

In summary the Mersenne and Sophie Germain primes are viewed by the operations that define them through the reduced residue system mod 18, 6 disjoint residue classes, LOPI 1, 5, 7, 11. 13, 17 mod 18 that contain all possible prime odd natural numbers (excluding 3, as already noted). By using the a broader equivalence relation derived from the mod 18 relation we show that the equivalence relations of mod 6 with LOPI 1, 5 represent the partitioning into two broader classes in order to define Sophie Prime or Non-Sophie Prime. The equivalence relation for each category avoids any overlap from these disjoint residue classes and therefore reduces the area in which we should look for Sophie Germain primes, as we know that if an integer meets the 5(mod 6) congruence relation then we know it could be a Sophie Germain, from LOPI classes 5,11,17. If it meets the congruence relation of 1(mod 6) we therefore know it cannot be a Sophie Germain prime as these elements all derive from LOPI residue classes 1,7,13, none of which can contain Sophie Germain primes. If we restrict our search for Sophie Germain primes to the residue classes 5, 11, 17 then we will know we have searched all odd natural numbers capable of meeting the requirements imposed by the operation. The equivalence relation 5(mod 6) is a better indicator of a possible Sophie Germain Prime from LOPI classes 5,11,17 mod 18, and 1(mod 6) an indicator of a non-Sophie Germain element LOPI classes 1,7,13, as well the LOPI/LEPI identity introduced at the beginning of this work can be used always for these indications:

LOPI/LEPI  N($_{a,bc..n}$)= (a + b + c...n).,= N≡ $LOPI/LEPI\ (mod\, 18)$  with Sophie Germain primes only in residue classes LOPI 5, LOPI 11 and LOPI 13 and Mersenne Primes only in residue classes LOPI 1 and LOPI 13, with LOPI 1 being generated from exponentiation of non Sophie Germain Primes, $2^{1+6(a)}$ and Mersenne LOPI 13 being generated from the exponentiation of Sophie Germain Primes $2^{5+6(a)}$.

**Perfect Number Analysis using the LOPI residue mod 18**

Given that all primes exist in the LOPI reduced residue classes of mod 18 we would like to analyse the Sophie Germain primes leading to Mersenne primes leading to perfect numbers **[22, 23]** in terms of the classes of integers that may work together to produce a perfect number:

$(2p -1)(2^{p-1})$

In the earlier discussion we showed that under the constant LOPI identity contained in each integer mod 18 that all Mersenne primes are represented only in LOPI class 13 and LOPI class 1 and that which class they exist within is determined by the type of prime from which they were made.

All Mersenne primes come from the LOPI congruence classes in two different groups. $2^{1,7,13}$ or $2^{1+6(a)}$, LOPI 1, 7, 13 all correspond to the class of residue classes in which we cannot find a Sophie Germain Prime, therefore we call these the non-S.G. group of residue classes, the classes congruent to 1(mod 6). As opposed to the Sophie Germain Prime bases for a Mersenne, $2^{5+6(a)}$, Primes in LOPI classes 5,11,17.

Depending on the source prime of the Mersenne prime, all Mersenne primes are in one of two LOPI congruence classes, LOPI 13 or LOPI 1:

$2^{1+6(a)}$    Sophie Germain Prime 5(mod 6)     Non S.G. Prime $2^{1+6(a)}$ 1(mod 6)

| | |
|---|---|
| $2^5$ - 1 = 31 N LOPI 13 | $2^1$ - 1 = 1   N  LOPI 1 |
| $2^{11}$ - 1 = 2047 N  LOPI 13 | $2^7$ - 1 = 127 N  LOPI 1 |
| $2^{17}$ - 1 = 131071 N  LOPI 13 | $2^{13}$ - 1 = 8191   N LOPI 1 |

These two classes contain all Mersenne primes possible, except for the one time case of prime 3, which yields the unique LOPI 7 element Mersenne prime,

$2^3 - 1 = 7$, all other $2^{3+6(a)}$ generate LOPI elements class 7 that are multiples of 7, as is in reality the Mersenne 7, but the result is one, so we have considered it a Mersenne prime originating from the power of 3, the unique prime in LOPI residue 3 mod 18, not coprime to 18.

This division in the Mersenne primes is important as it tell us something about the perfect number that will be formed based on LOPI combinations under the definition for a perfect number. The congruence 3(mod 4) however does not distinguish between these two sets of integers that may become Mersenne primes, as shown earlier in this work. It finds a congruence relation for all the LOPI classes, 1,5,7 and 11, while does not work for integers in the LOPI 13 and 17 classes. This is a problem that has shielded more partitioning work in terms of Mersenne numbers and their relationship to Sophie Germain primes and the resulting perfect number from the combinations of a Mersenne NON-SG in the exp.-1, and a Mersenne $^{SG}$-1. A separate congruence relation can be used to distinguish between the two exponent types as 1(mod 6) Non SG exponent primes, and 5(mod 6) for exponent SG M primes.

This becomes important in understanding the combinatorial relationship as the root prime of the Mersenne is indicative of what happens to the Perfect number upon +1 and if the Perfect number LEPI 10 is a possible Euler totient function, which may also indicate to us that by adding one to Perfect numbers made from the two variations of Mersenne primes.

Mersenne primes exist in the LOPI congruence classes 1 and 13, LOPI 1 being generated by the exponent elements of $2^{LOPI\ 1\ +\ 6(a)}$, LOPI classes 1,5,11, LOPI N congruent to 1(mod 6).

LOPI 13 Mersennes being generated by $2^{LOPI\ 5\ +\ 6(a)}$ from the LOPI which are classes that can make a Sophie Germain prime, or from the LOPI classes 5, 11, 17, LOPI N congruent to 5(mod 6).

Congruence relation 3(mod 4) finds true for LOPI 7 as well, but cannot be a Sophie Prime. As well we see it holds for 19, a LOPI 1 element that also cannot be a Sophie Prime, but it is a Non Sophie Germain Mersenne. Here we therefore divide the Mersenne primes based on the equivalence relations of Sophie Germain prime is congruent to 5(mod 6).

Non Sophie Germain prime is congruent to 1(mod 6).

We will analyse the components of the perfect number equation to show that all perfect numbers will be in LEPI 10.

First we apply the operation using all the LOPI residues to see the possible LEPI classes into which the results will fall.

$2^{1-1} = 2^0$ = LEPI 10    $2^{5-1} = 2^4$ = LEPI 16    $2^{7-1} = 2^6$ = LEPI 10
$2^{11-1} = 2^{10}$ = LEPI 16    $2^{13-1} = 2^{12}$ = LEPI 10    $2^{17-1} = 2^{16}$ = LEPI 16

Here we can see that each LOPI residue when carrying out the operation of $2^{p-1}$

yields LEPI 10 for exponents belonging to the Non Sophie Germain LOPI classes, $2^{1+6(a)-1}$ and elements in LEPI 16

for exponents belonging to the Sophie Germain LOPI classes, $2^{5+6(a)-1}$.

Now we are down to two LOPI classes of elements in which the Mersenne primes exist:

Mersenne $2^{S.G.}$, $2^{5+6(a)} -1$ = LOPI 13   and Mersenne $2^{Non\ S.G.,\ 1+6(a)}$ = LOPI 1

 S.G. exponents $2^{5+6(a)-1}$ = LEPI 16    and Non S.G. exponents $2^{1+6(a)-1}$ = LEPI 10

If we put these into their linear progressions:

13 + 18(a)(16 + 18(b) = 208    LEPI 10      a≥0
13 + 18(a)(10 + 18(b)=130 LEPI 4
1 + 18(a)(10+ 18(b)=10 LEPI 10
1 + 18(a)(16 + 18(b)=16 LEPI 16

When we put these four combinations together to form the perfect number we know that in $(2p -1)(2^{p-1})$  p
 must equal p in both parts.  From this we can see that the two linear binomial progressions will not work are:

13 + 18(a)(10 + 18(b)=130 LEPI 4
1 + 18(a)(16 + 18(b)=16 LEPI 16

or            $2^p-1$                              $(2^{p-1})$

Mersenne ($2^{S.G.,}$ or $2^{5+6(a)} -1$ = LOPI 13) times ($2^{1+6(a)-1}$ = LEPI 10),

here  the exponent prime $2^{5+6(a)}$ cannot equal the exponent prime in $2^{1+6(a)-1}$ as they are from different LOPI

linear progressions, and one is Sophie Germain possible LOPI classes, 5,11,17 and the other part is generated from the other twin set of LOPI congruence classes, non Sophie Germain containing LOPI families, 1,7,13.

Therefore the primes never equal each other.

 The combination prime in the Mersenne does not equal the prime in the other part of the equation, but we get a Mersenne prime in the LOPI classes 5,11,17 multiplied by a prime from the other twin classes of these or the non-Sophie Germain classes, $2^{1+6(a)-1}$, LOPI classes 1,7,13, therefore these primes will never match up and cannot meet the requirements of the operation to produce a perfect number. The LOPI residue classes of mod 18 in the exponent for the Perfect number definition are generated every 6 in the exponents of the operation, therefore this base combinatorial example carries through to the elements in all the LOPI classes as they are grouped by the residue and all elements contain the residue that controls what class of element will be generated by these operations that form a perfect number.

The second case, as well, does not work as well as the prime exponents do not equal each other as they are the partitions of Non Sophie Germain LOPI classes and Sophie Germain Classes:

(1 + 18(a))(16 + 18(b))=16 LEPI 16

$2^p-1$                             $(2^{p-1})$

Mersenne $2^{\text{Non S.G.}}$ or $2^{1+6(a)} -1$ = LOPI 1) times ($2^{5+6(a)-1}$ =LEPI 16)

We have the Non Sophie Germain LOPI residue classes forming the Mersenne and the Sophie LOPI residue classes forming the latter part of the equation, therefore these will not produce a perfect number under the requirement that p=p in $(2^p-1)(2^{p-1})$.

Now to list the two that will work:

(1 + 18(a))(10+ 18(b))=10 LEPI 10   Poly 10 + 18(b) + 180(a) + 324(ab)
(13 + 18(a))(16 + 18(b)) = 208      Poly 208 + 234(b) + 288(a) + 324(ab)

Only these two possible linear combinations yield the necessary requirements of Perfect = $(2^p-1)(2^{p-1})$ where the prime in the Mersenne equals the prime in the other part of the equation.

Mersenne $2^{\text{S.G.}}$, $2^{5+6(a)} -1$ = LOPI 13   and $2^{\text{S.G. }5+6(a)-1}$=LEPI 16
Mersenne $2^{\text{Non S.G.}, 1+6(a)}$ = LOPI 1    and $2^{\text{Non S.G. }1+6(a)-1}$=LEPI 10

(1 + 18(a)(10+ 18(b)=10 LEPI 10  Non-Sophie Germain type Perfect
(13 + 18(a))(16 + 18(b) = 208  LEPI 10  Sophie Germain type Perfect

We therefore put forth the Theorem based on the analysis of the LOPI residue classes mod 18:

All perfect numbers are in LEPI class 10 derived from the combination of:

(Mersenne $2^{\text{S.G.}}$,-1) with $2^{5+6(a)} -1$ = LOPI 13   and $2^{\text{S.G. }5+6(a)-1}$=LEPI 16

(Mersenne 2 $^{\text{Non S.G.,}}$ 1) with non $2^{1+6(a)}$ = LOPI 1   and 2 $^{\text{Non S.G. }1+6(a)-1.}$=LEPI 10

(1 + 18(a)(10+ 18(b)=10 LEPI 10     or    (13 + 18(a)(16 + 18(b) = 208  LEPI 10

For Perfects of the base Non-Sophie Germain we conjecture that upon + 1, these perfects are Euler totient values for the next prime 11 LOPI element.

Perfect Numbers base Non-Sophie Germain + 1= Prime LOPI 11

And Perfect number base Sophie Germain + 1 meet equivalence 0(mod 7), are multiplies of 7.

**Twin Primes:**

Twin primes, p + 2,  can be investigated completely within the LOPI 17 residue class mod 18, as all subgroups within the Superset 17 are made by the twin LOPI multiplication seed values, 17, sub 17*1 with b>0, 35, and 143.  Given the fact that these convolutions will contain the product of a twin prime pair all within LOPI 17, we already have a huge advantage in partitioning based on the LOPI Mod 18 system.  Once we perform the convolutions we know that the elements within the subsets are the products of the twin primes.  Also we can construct the product by letting a=b for subsets based on LOPI 7*LOPI 5, and LOPI 13*LOPI 11 thereby setting up the condition that the two are twin elements and perhaps are a prime elements in twin LOPI congruence classes. Also for LOPI 17 * LOPI 1, since this is the pair within mod 18 that leads to the twins LOPI 17-LOPI 19 (LOPI 1 +18), within this convolution we let a=b+1, since they are a cycle off in their prime twin elements. Twin primes exist within mod 18, in LOPI sets [5,7], [11,13], [17,1].  LOPI 17 elements' twin is not in mod 18, but the next cycle, mod(36). We therefore construct the convolution to yield all results with a=b+1, (1+18(4)(17+18(3)=5183, as expected LOPI 17, returns a multiplicity of one from the query of 5183 in the range needed of the subset (17*1) LOPI non-prime subset. When we construct these polynomials we fix the relationship to be such that a and b will be elements in a non-prime subset of LOPI 17  and the product of this convolution will yield a multiplicity of 1 in the subsets if it is a prime element times a prime element, therefore a twin prime times a twin prime. Here we can see that since *a* is an element of the natural numbers, N, and b as well, even when it is defined as a+1, we see that in these results we will always come upon a set of twin primes in our convolution when we see a multiplicity of 1 in the particular LOPI subset for which we constructed the product to be generated. Primes occur in the absence, the misses, where the seed operation "gears" when treated as convoluted binomials cannot produce a factor in the cycle of adding 18 that produces a prime. The only operation that yields the prime is the original Superset LOPI progression produced by the LOPI +18(a), with a≥0, not to be confused with the non-prime subset with the LOPI and LOPI 1 subset in convolution as discussed above. The primes do not undergo a new function within the Superset A, but remain unchanged, while all the non-primes undergo an equivalence relation as discussed above to be partitioned into their LOPI non-prime subsets.

The fact that the pattern for the calculation of the factors of the non-primes exists and once calculated by the convolution discussed above, we are only limited by the power of our CPUs and the coding of the attached algorithm for twin prime query, to know the value of large secondary primes and twin secondary primes and their factors. Also given that the LOPI non-prime subsets can be subtracted out from the LOPI super-sets to reveal the pure odd prime subset for each LOPI congruence class, therefore the complete prime set in the Universal set of all the collective Super-sets can be known.

A pattern of factor generation producing the Non-Prime elements and skipping the Prime elements in N is revealed as being inside the 6 LOPI linear binomials created by using the reduced residue system Modulo 18, the totatives of ω(18) and by the further convolution of specific Dirichlet linear binomials to generate the factors of each non-prime within each of these 6 LOPI congruence classes. The primes and non-primes add rigidity and flexibility to structures in nature, much as we have different strengths and types of bonds in Chemistry, some that allow for long distance coordinating weaker interactions where the distance is not correct to form very stable covalent bond, the numbers begin appear in this manner as well. The primes multiply to give us a little less rigid structure with which to work, as earlier Euclid called the plane elements, once two linear primes are multiplied by another prime, we get a number that does contain the ability to be divided into two wholes, a secondary prime. As we continue through the mechanism through convolution, we generate all possible factors leading to the non-primes, the "non-linear", the plane and solid integers, as described in book 7 of Euclid's The Elements. We see that all the numbers contain a lowest odd or even partition integer by which we can create these families, we can partition these elements into disjoint Super-sets, that guided by their LOPI and mod 18 predict certain behavior from the basic interaction of this core identity. The primes, being the gaps of the composite forming polynomials, allows us a structure that we can't really change from which to adapt, grow and change direction linearly. In generating the non-primes through convolution, however, LOPI modulo 18 gives us a growth model for primes and non-primes, a way that shows the composites are made by the linear binomials from the LOPI residues mod 18, and the primes just aren't constructed from factors generated since we do not allow their one factor to be used, 1. The primes are therefore not mapped over into the non-prime LOPI subsets. The linear binomial progression LOPI +18(a), creates primes and non/primes and from the lowest factors of the seed non-primes, all the non-prime elements of the LOPI congruence classes grow, while the primes continue only in their linear progression toward infinity, with the convolution equivalence relation not being able to change them within their Super-sets. The Prime subset in A is the complement of the total composite subsets which we calculate in A with the LOPI residue system mod 18, or Superset A\LOPI subsets=S'=Prime subset in A, therefore the possibility is zero that a prime can be in the subset since it is an element only in the complement of those LOPI subsets of the LOPI Superset congruence class. Within this Modulus system an integer can easily be accessed as to its primality and number of divisors and identity of its factors as shown above, once the convolution is complete for each LOPI congruence class. Even without completing the final convolution, the equations can be used independently to create large non-primes with known factors and primality, whether by the polynomial or the convolution. Within this system of polynomials perhaps also the Reimann hypothesis can also be analyzed at a later point in order to understand how the functions may work together to reflect the Reimann zeta function in correlation with the gaps generated in the non-prime polynomials presented through the LOPI mod 18 system.

The polynomials generate only the composite numbers in the set of integers, N coprime to 18, therefore the gap in the codomain of the range are the primes being cyclically skipped.

Seed Constant + LOPI*18(b) + LOPI(18(a) + 324(ab)≠LOPI + 18(c)≠a Prime LOPI element

Seed Constant + LOPI*18(b) + LOPI(18(a) + 324(ab)=LOPI + 18(c)=Non-prime LOPI element

Twin primes exist within mod 18, only in LOPI class pairs [5,7], [11,13], 17's twin is not in mod 18, but exists in the next cycle of 18. In fact LOPI 7 can either be prime or non-prime within the cycle of 18, and also independently LOPI 5 can be either prime or non-prime, to say that these two congruence classes cannot possibly be prime at the same time now or in the infinite future is to imply that even at the low values of LOPI 5 and LOPI 7 they cannot both be prime, which we know to be untrue. Therefore we can extend this basic relationship to infinity because if not then we must reverse the logic and say that a number and a number +2 changes the degrees of freedom of the other and if this is true, it will always be true, even in the other direction, so if LOPI 5 + 18(a), where *a* goes to infinity being ≥0, cannot be prime at the same time as LOPI 7 +18(a) where *a* goes to infinity by ≥0, then LOPI 5 + 18(a) can never be prime when LOPI 7 +18(a) is prime, even at a=0. We know this to be false and now with the work above we see that "primality" in a LOPI congruence class is dependent upon the mechanism of factor generation. Also, the gaps created by the LOPI Mod 18 system during convolution is not dependent upon their twin neighbor, but upon the iteration of 18 from each LOPIs lowest non-prime product factors which multiply to contain each LOPI core identity. We would have to say that non-prime generation for different LOPI subsets would be equal in the way that they would always co-produce non-primes, or always co-produce primes, when in reality we show that these congruence classes by the binomial convolution produce different results, always we can say, a P*P, P*NP, NP*P, NP*NP for each subset LOPI 17, which are 3 sets of two linear binomials convoluted. If we contend that this pattern will change then we will not produce all the non-primes since we say that we can't have a P*P in the future, only P*NP, NP*NP,NP*P, but this isn't possible as we prove that elements that are secondary primes in the subsets grow from the primes, the tertiary grow from the secondary, if we remove the prime*prime we remove the P*NP, which removes the NP*NP in the pattern of factor generation.  This doesn't exist, as the secondary prime subsets exist to infinity since it has been proven that primes exist to infinity and therefore new secondary primes will always exist by the new prime created and the past ones already generated in the LOPI Dirichlet linear progressions mod 18. Therefore we argue that as a partitioned secondary prime subset will always be populated to infinity that represents the secondary primes, then so too will the twin prime secondary primes as shown to already exist in the cycle of LOPI convolution for families [5,7], [11,13], and [17,1]. The fact that they are twin is incidental within the pattern of each LOPI congruence class, each prime element location within each LOPI family is determined by the convolution pattern of non-primes within each LOPI and at the jump where 1 was needed as a factor to generate the prime. The fact that this Prime/Non-prime pattern is being created from the beginning of the first non-primes and further non-prime elements in each LOPI congruence class are being generated by a cyclic degree one polynomial with two natural number variables, a and b, implies that this cyclic pattern will continue changing based on the growth of variable a and b, natural numbers ≥0.  Important, also, is the fact that as this cyclic pattern of non-prime generation is occurring, the cyclic skipping of the primes is also generating a pattern, the pattern in which the primes exist in relation to their non-prime family elements.

If we c set of each of the subsets LOPI 17, by which a=b, and for 17*1, a=b+1, if Prime*Prime then these will be the elements of each subset that appear one time as secondary primes do, we know of these generated elements, twin P * twin P, when we set a=b, that this product element in order to be twin prime factors should appear only one time in the subset operation of this LOPI 17, if it appears in another of its LOPI subset or again later in the same subset, then the multiplicity is not one and  we know we do not have twin P* twin P product. For example, if we look at data matrix where column = row of subset generated by 7+18(a)(5+18(b), I get possible candidates for twin P*twin P, I can cross reference with a query to see if they are in the other two subsets of LOPI 17, and then I know the primality of this element from my subset of LOPI 17. The fact that all possible combinations for twin primes result in the LOPI congruence class 17 non-prime subsets, is amazingly convenient for potential twin prime generation and multiplicity/factor search.

It's interesting to note that although the number 3 is considered prime, as a LOPI to generate primes cannot be used as it generates more and more values of LOPI 3 continuously upon LOPI +18(a). But LOPIs 3, 6 and 9 generate some Sophie primes but never a Prime within their LOPI, perhaps the number 3 needs its own new definition. It needs more study as a unique component of mod 18, $2*3^2$, where $3^2$, cyclic group order 6. The continued use of three as a basis for looking at the continuation of progressions of primes maybe creates confusion in pattern recognition. Here with modulus 18 we show a clear pattern of partitions and sub-partitions, groups and subgroups, that agree nicely with the past work on twin primes, Sophie Germain primes and the Prime Number Theorem, and furthers theory to reveal a new pattern in the generation of all non-primes with the totatives of ω(18), or the LOPI mod(18), which reveals the pattern of the primes as well. They are a subset under each LOPI, but are not formed by an equivalence relation, R, on A, by the convolution of particular linear Dirichlet progressions. We therefore have all the non-prime element subsets partitioned from the larger partition LOPI residue identity and once partitioned by the equivalence relation the primes remain as the last subset in the Superset.

{$LOPI_1$ Superset} – {$LOPI_1$ Subsets from Convolution}= {$LOPI_1$ prime subset}

The pattern of non-prime generation within the odd LOPIs that can generate prime numbers follows a cyclic pattern, one distinguished by the Lowest Odd Partition residue identity of each group. A large number can be classified immediately into its LOPI congruence class without knowing anything else about it by the Identity LOPI N≡LOPI (modulus 18), or we can add all the summands together and get LOPI as well, stopping at LOPI/LEPI depending on the odd/even nature of our number. We don't need to know any more than the LOPI to open a door to an amazing amount of information about that number, as shown above. We have many ways to use and generate the information revealed by a numbers LOPI/LEPI. The binomial convolutions and the polynomial can all be treated and investigated to reveal the nature of the range and possibilities of primality, twin primes, secondary primes and Sophie Germain primes. The LOPI reduced residue system mod 18, and the complete covering system mod 18, reveals a very elegant system that goes in harmony with the known theories and lends further support to them, opening new partitioning information and providing a structure from which to investigate further various aspects of our number structure. Non-prime subset multiplicity and the mechanism of product formation after factor generation using the convolution of LOPI Dirichlet linear progressions give us a deterministic tool to evaluate the image/anti-image the primes and non-primes share in each congruence class after further partitioning.  We can create the algorithms attached to this article to go further with confidence that we are going to find what we are looking for, as with the twin discussion whose products are within subset 17, if we restrict a=b  and for a=b+1 for LOPI 17 NP subset ((17 + 18(a))(1+18(b))) in the subsets of LOPI 17 and then query the three subset matrices in the grid range discussed, an output of multiplicity of one, then we know we have just found a twin pair secondary prime product. Also by the other type of analysis generating by the Delta, where the factors are less readily available, we can generate directly the list to query for multiplicity one in the subsets of LOPI 17. The delta generation equation for the non-prime elements where both a and b are iterated each time for subset 7*5, and 11*13 with a=b, is $N_{x-1} + 18(a_{x-1} + 36) = N_x$ and for the potential twin product non-prime elements in subset are generated by  17*1 past the second iteration product, 1295, to take into account the imbalance in the cycle, is  $N_{x-1} + 18(a_{x-1} + 36)$, from 1295 we have 1295 + 18(36+54)=2915 and continue as normal, noting the first products, 35 and 323, where 323 is an already known twin prime product, 19*17, with multiplicity one as expected in LOPI 17 non-prime subset (17+18(a))(1+18(b)) with b>0 and for possible twin prime product production, a=b+1.

Later also with more analysis and modeling, such as Fourier transform modeling we can see exactly how the convolutions of each LOPI are interacting, with each other and independently, maybe there is a pattern to these 6 congruence classes, and maybe a larger pattern to the 9 and then to the complete covering system mod 18.   If there is an overarching interaction or a symmetry/asymmetry component to the total system mod 18, perhaps this system allows us to have symmetry and symmetry breaking in order to form new adaptations and structures, and this distribution of primes and non-primes is a perfect balance of flexibility vs. rigidity at a finer level that promotes change and stability. It would seem so, as Euler stated nature seems to have divided the primes into two classes, 4n + 1,  4n-1, **[24, 25]**, by using a partitioning modulo to exclude 3 and 2, we see the primes are divided evenly among 6 congruence classes, perhaps this pattern reflects an effect on the possible outcomes.  For instance maybe when two quantities are just right, or wrong, their product stops short or goes too long or opens a door to another reaction, that if the resulting structure had contained no lines of symmetry for division by integers, perhaps it would have stopped the reaction from continuing or triggered another response from the surrounding system . Our biological world appears to have pattern underpinnings in the structures possible by different types of integers, as we see here via the LOPI /LEPI mod 18 they are cyclically produced. **[26]**

At any rate, the natural numbers are partitioned into 9 odd classes and 9 even classes by mod 18, and further the odd classes are partitioned into 6 classes of mixed prime and non-prime odd integers by the reduced residue system mod 18, which are represented by 1(mod 6) for LOPI 5,11,17, possible Sophie Germain primes, and 5(mod 6) for LOPI 1, 7, 13, square containing odd LOPI classes. Each set of classes contain one of the twin prime classes and each set containing only one possible Sophie Germain Prime, as discussed above, offer a wealth of information when partitioned further into their unique LOPI reduced residue classes mod 18, {1, 5, 7, 11, 13, 17}.

Primes occur in the absence, the misses, where the seed operation "gears" when treated as convoluted binomials cannot produce a factor in the cycle of adding 18 that produces a prime. The only operation that yields the prime is the original Superset LOPI progression produced by the LOPI +18(a), with a≥0, not to be confused with the non-prime subset with the LOPI and LOPI 1 subset in convolution as discussed above. The primes do not undergo a new function within the Superset A, but remain unchanged, while all the non-primes undergo an equivalence relation as discussed above.

The fact that the pattern for the calculation of the factors of the non-primes exists and once calculated by the convolution discussed above, we are only limited by the power of our CPUs to know the value of large secondary primes and their factors upon query in the known grid area of the known subset data matrix . Also given that the LOPI subsets can be subtracted out from the LOPI Super-sets to reveal the pure odd prime subset for each LOPI grouping, therefore the complete prime set in the Universal set can be known. By recognizing the factor generation pattern based on the core LOPI that all P/NP contain, the Non-Prime elements and the Prime elements are revealed as being inside the 6 LOPI Dirichlet linear binomial progressions created by using the reduced residue system Modulo 18, the totatives of ϖ(18), and by the further seed convolution of specific Dirichlet linear binomials to generate the factors of each non-prime within each of these 6 LOPI congruence classes. We see as shown above,

0≠ Prime LOPI elements in (LOPI + 18(c)) — seed + LOPI*18(b) + LOPI(18(a) + 324(ab)

0= NP LOPI elements in (LOPI + 18(c)) — seed + LOPI*18(b) + LOPI(18(a) + 324(ab)

and for prime the equation equals the original grouping linear Dirichlet linear equation for the LOPI, LOPI + 18(c).

In the following section we wish to analyze the well known classic Goldbach Conjecture using the LOPI reduced residue mod 18 system to introduce templates showing that prime solutions will exist in the Dirichlet Linear Progressions needed to create each even natural element ≥2, not using 2 or 3 as prime, as they are not LOPI reduced residues mod 18. We can say that LOPI $_1$ + 18(a) + LOPI $_2$ + 18(b) = LEPI$_3$ + 18(c), the linear equations being those of the Dirichlet LOPI linear progressions that populate each LOPI congruence class mod 18 with primes and non-primes to infinity will provide a template for all even natural number. LEPI 2, LOPI 1 + LOPI 1 is not prime, but still a very needed LOPI template mod 18 that leads to the linear progression containing needed primes for larger LEPI 2 integers. For LEPI 4, 3+1 template will not be used as LOPI 3 does not generate primes upon iteration, therefore we pass to (11 + 18(a) + (11+18(b)) template to generate LEPI 2 elements >2. LEPI 6, template seed (1+18(a))+(5+18(b)) is needed and again we do not use 3 as a LOPI, it is not a coprime residue to mod 18 and will produce no other primes than itself in the LOPI congruence class 3, mod 18. We rewrite the strong Goldbach Conjecture in terms of the Dirichlet LOPI Linear Progressions mod 18, as:

All Even elements in the natural numbers can be written as the addition of two primes within the reduced residue LOPI linear progressions by LOPI P + 18(a) + LOPI P + 18(b), with least element even residues 2, 4 and 6 not made by two coprime residues to 18. Upon multiplication of the LOPI elements produced, as shown below, the product can be queried in specific locations of the known subset data matrix to confirm a multiplicity of one, thereby confirming a secondary prime product and a Goldbach prime + prime solution.

**Goldbach Conjecture**:

The above system will now be used to provide congruence templates based on the basic arithmetic relationships expressed as the LOPIs in the linear binomial form. Templates are shown to produce all even numbers as conjectured by Christian Goldbach in a letter to Euler in 1752 **[27]** suggesting that all even numbers >2 can be written as the sum of two prime numbers, known as Goldbach's Conjecture . In working with all the information revealed by LOPI/LEPI reduced residue system mod 18 templates can be formed to yield congruence relationships to obtain all the possible LOPI + LOPI = LEPI solutions, combining this system with the Prime Number Theorem and the Dirichlet Theorem of infinite primes in reduced residue systems and by taking advantage of the core LOPI identity to create the LOPI templates needed in the residue classes, we will proof for the Goldbach Conjecture along with experimental data showing that the primes are available and sufficient to yield pairs of LOPI P/NP + LOPI P/NP = LEPI N. They do exist and can be generated and proven to be P + P by the matrix query for a multiplicity of one in the particular subsets of the respective LOPI. Also to attempt to show that it is impossible that a Prime + Prime=LEPI not occur, due to the mechanism by which factors are generated and non-prime products formed. Finding LOPI/LEPI, a mathematical identity coinciding with the six congruence classes of the reduced residue system modulus 18, yields an elegant system in factor generation and the prime/non-prime relationship. According to Bertrand's Postulate, [28] we also have proof that at least a prime will be present in each linear binomial, as it is n<p<2n, putting limits on the non-prime consecutive elements.

(LOPI + 18(a)) + (LOPI + 18(b)) = (LEPI + 18(c)), with c or c-1=b to get the maximum number of cycles of 18 needed to create all the possible P + P pairs.

These equations will give the whole set of pairs that exist in the templates of the $LOPI_1$ + $LOPI_2$ = LEPI, which are the Dirichlet linear progressions of modulo 18 described above.

The needed linear binomials can be known based upon the basic mathematical rules of addition that the core LOPI follow throughout the LOPI congruence classes.

The lowest even partition identities of modulo 18:

0   2   4   6   8   10   12   14   16

18   20   22   24   26   28   30   32   34

Upon iteration of +18 we fill the congruence classes as above with the reduced residue system of the LOPI, but here it is not the reduced residue system, but still within each class we do have congruency.

Now we will list them as their linear binomial:

0 + 18(a)  2+18(a)   4+18(a)  6 +18(a)   8 + 18(a)   10+18(a)  12+18(a)  14+18(a)  16+18(a)

These linear binomials will populate each LEPI congruence class to infinity, as a is the set of N, from 0 to infinity.

Now we see exactly how we reach these values by 1, 5, 7, 11, 13, 17, the reduced residue classes of modulus 18 that contain all possible odd primes in the natural number set, as shown before in this research. Also we know that the non-prime generation of elements within each congruence set is evenly distributed and cannot create consecutive non-primes in a whole cycle of 18, all 6 numbers of the reduced residue classes on an interval of 0-18 cannot be all non-prime, as we have shown, therefore we know that on any given interval of 0-18(a), we will have at least one prime per congruence class less than the even number we are inquiring about. We know that within each of the congruence classes we have an even distribution of primes according to the prime number theorem and we can use the formula also to give us a rough idea of how many primes we have to meet the Goldbach Conjecture.

The LOPI list of combinatorial relationships that yield these particular LEPI congruence classes upon (LOPI + 18(a) + LOPI + 18(b))=LEPI + 18(c).  For example,

LEPI 2, we use LEPI 20, since 2 is a LEPI residues classes mod 18

But we will need to use 1, since it is our LOPI reduced residue Mod 18.

1+1=2

17+3=20, but here is a 3 and in the LOPI/LEPI mod 18 system only reduced residue system mod 18 is used as a LOPI 3, will never yield a prime, other than itself, that has been categorized a prime, but for us we extend the relationship of the core LOPI to infinity, and therefore LOPI 3 will always yield a non-prime. Therefore the template 17+3 is not valid, nor 5 + 15.

7+13=20, a P + P, a valid template.

11+9=20, again no, we have a LOPI 9, a non-reduced residue element of mod 18, so again a no, it will never produce primes we can use for the Goldbach Conjecture.

For LEPI 2, we have two templates of two different linear binomial Dirichlet progressions which guarantee some evenly distributed primes to infinitude, by which to extend into space by iteration of 18 in each LOPI which will create all even numbers in the LEPI Congruence Class.

1+1=2,   7+13=20  therefore can be written as :

$2 + 18(c) \equiv 18(a) + 1 \pmod{1+18(b)}$ , remembering that in the convolution 1 is never allowed.

$20 + 18(c) \equiv 18(a) + 7 \pmod{13+18(b)}$

We can also write these as before mentioned:

$(13 + 18(a)) + (7 + 18(b)) = 20 + 18(c)$.

$(1+ 18(a)) + (1+ 18(b)) = \text{LEPI} + 18(c)$.  Later I will change the variables to reflect that use of the c in the LEPI.

For the congruence relations above lets continue.

We start by getting the value of *c* for our even number,

92, 92-2/18=5, we therefore know that 92 will have a LOPI prime + LOPI Prime in the LOPI reduced system of mod 18. We can get all those pairs by starting from a=0 and b=c-1 for all LEPI<9, therefore we have and a goes up by 1 and b continue to go down by one:

$2 + 18(5) \equiv (7 + 18(0)) \pmod{18(4)+ 13}$, = 92=7+85

To continue we get

92=7+85

92=25+67

92=43+49

92=61+31

92=79+13

We go until b=0 and we started with a=0 in the linear binomials.  With just this one template at such a low number as 92, we have two LOPI Prime + LOPI Prime=LEPI.  From the simple identity of the LEPI we know the number of different solution pairs and we know how to get them. For LEPI 2, element 92-2/18=5, the five pair solutions above, two of which out of five meet the Goldbach Conjecture requirement, Prime + Prime = Even LEPI. As numbers grow we get an even larger see-saw effect balancing between a=0, in the residue LOPI linear binomial all the way until b=0 in the modulo Linear binomial. The range becomes larger and larger for each congruence class the larger N becomes, thereby giving us even more opportunities to obtain prime pairs to meet the Goldbach Relationship. These congruence classes come from the LOPI of the reduced residue system Modulo 18, the primes are present by the Prime Number Theorem and due to the fact of the factor generation mechanism, which creates gaps in the range of the codomain, the primes are those gaps in the process of creating non-primes by the convolution of the specific linear binomials discussed above.  The non-primes subsets are equal in number of elements for each LOPI as they each contain 3 subgroups of Non-primes, the LOPI with four, are really 2 full subsets and 2 twin ½ subsets when the $LOPI_1=LOPI_2$.  Each non-prime subset is generated by the same rate in each cycle of the linear convolution, the only difference is the residue which is still between 1-17, therefore for cycles of 18 the rates of non-prime generation are the same within each subset equal to $(LOPI_1 + 18(a))(LOPI_2 + 18(b))$, with a,b≥0.  The LOPI mod 18 system provides evidence that the primes are distributed evenly throughout the natural odd numbers since the non-primes are generated by the same mechanism of iteration and convolution of 18 no matter the LOPI congruence class. They run parallel to each other and only differ in their point of origin from zero, in a cycles of 18, of (0)+1, 4,2,4,2,4,+1+(18) to complete the loop 0-18. Each LOPI will have roughly an even distribution of primes amongst them,  as they are generating non-primes by the same mechanism with the same modulo 18 iterated as discussed in detail below.

All LOPI elements generated by the Dirichlet linear progressions LOPI + 18(a) and all non-primes generated by (LOPI$_1$ + 18(a))(LOPI$_2$+18(b)), a,b≥0 for the general case.

Below are the templates for the Goldbach/Dirichlet see-saw partition congruences made from the convolution of the linear equations with even and infinite distribution of primes, as they are running at the same interval and rate we know their number of elements in each subset are equal in the cycles of Infinity(18), since no matter how many cycles we go, we produce 6 elements, 1 each for every cycle of 18, as well for non-primes we are creating equally in each subset nxm, or a+1(b+1) number of elements, with the subset for the squares, where a=b, ½ of nxm, therefore we are regularly producing primes as we generate the non-prime elements in the subsets where the factors are never one and cannot multiply in a 2n+1 manner, but in an 18(n)+ LOPI, in a sequential manner but much greater, + 18(5), 18(7), 18(11), 18(13), 18(17), showing that the primes are present, cyclically and consistently, the gaps get bigger because the cycle of 18 get bigger and the number of non-primes continues to grow more quickly than the rate of the prime numbers.

Integer 13484, LEPI 2, that is within the scope of what we have calculated 13484, we know their will be Xa number of elements in each Goldbach/Dirichlet See-Saw Partition, for 13484 or in 0-13482, Xa=18(749) + LOPI or up to cycle 13482, 18(349), c=749

(2+ 18(749))≡(1+ 18(1)) (Mod 1 + 18(747) a=1, b=c-2, for LOPI 1

(2+ 18(749))≡(7+ 18(0)) (Mod 13 + 18(748) = 13,484=7+13477

We will have 748 solutions in the linear Dirichlet equations for LOPI 7, and 374 possible solutions with LOPI + LOPI = LEPI 2. The Goldbach/Dirichlet See-saw Congruences go from a=0 to b=c, c-1,

The first pair then is 7, 13477 and in the data table, both confirmed prime, as absent in any subset of LOPI 13. Another pair confirmed easily is from the template LOPI 1 + LOPI 1=LEPI 2, with (1+18(4) + (1+ 18(745))=13484. Here we already have two from the LOPI Modulus 18 templates. The gap between a=0 and b=c-1, reveals the number of pairs we will have to choose from. These are the disjoint reduced residue LOPI Supersets of LOPI 13 and LOPI 7 and primes distributed within these reduced residue congruence classes of Mod 18. Below is the long way to show the formula that produces the consecutive non-primes for each subset and type of iteration that is in line with the non-prime/prime pattern generation theory via linear binomial convolution LOPI mod 18, as the gears of non-prime production leaves gaps upon iteration of 18.

All the LOPI element non-primes are generated from a few seed non-primes in each LOPI reduced residue class. From these seeds, we make the non-primes by iteration of 18 onto each, one or the other and then multiply. This gives us all the non-prime elements in each of the residue classes of Modulus 18. By choosing mod 18 and the residues, which coincide with the LOPI identity, we assure that no odds from LOPI 3, 15, 9 can be included and also we ensure that no integer from these LOPIs are ever generated, as well 1 as a factor is prohibited as to avoid any primes appearing in subset generation for non-primes with non-trivial factors.

In order to offer a numeric proof, we can calculate the number of non-primes generated by each of our choices from our seed operations using the following method.

LOPI 13 has four seed operations as discussed earlier:

7 x 7 = LOPI 13

13 x 1 = LOPI 13, with b never 0

11 x 11 = LOPI 13

17 x 5 = 85

For each of these we add 18 to both of the seed factors, to the right and to the left as such, knowing that the highest integer for that LOPI subset can be set by the highest chosen a and b for iteration of +18 to both, such as in the example below, we chose a,b=1 for both to look below 805 with all the subset seed generators.

(17+18)(5+18)=805

(17+18)(5) =175, with every non-prime in this line will equal previous + 18(5), or $175_{a1}$+90=265.

(17)(5+18) =23(17)=391, with every non-prime after that equal to 18(17)+ $a_{n-1}$, 391+306=697.

We can know how many non-primes exist based on these seed operations and their iterations via +18 in the prescribed way above.

We get a delta for each of these infinite lines of non-prime elements for each subset.

Delta both=$18(36+a_{n-1)})+ N_{x-1}$, the first a=$X_1$-$X_0$/18   (17+18(1))(5+18(1))=805 our chosen upper limit

Delta Right iteration, or 0, 1=18(element unchanged) or 18(left unchanged) + $N_{x-1}$, or $N_{x-1}$ + 306=$N_x$

Delta Left, or 1, 0 = 18(right unchanged) + $N_{x-1}$ or 90 + $N_{x-1}$

Here we see that from the left LOPI iterated, 17 and the right 5 not, our Delta is 18(5)=90 and up till 805, yields:

(17+18)(5)=175, for LOPI 13, + 90 every time

85+90=**175**+90=**265**+90=**355**+90=**445**+90=**535**+90=**625**+90=**715**+90=**805** or 8 np elements LOPI 13

For the 0,1, of (17+18)(5+18) with Left 17 held steady and right 5 not, we get Delta 18(17)=306 each time;

**85**+306=**391**+306=**697**+306=1003, we get 2 numbers only for LOPI 13 in this seed on this arm.

In total we see 10 non-prime elements of LOPI 13 are created from the seed 17*5 seed up to 805.

We do this for each seed operation to get the total non-primes in order to subtract them from the total elements in LOPI 13 up to 805 as a check and proof that we generate all non-prime elements in LOPI 13, not including any 3, 9, 15.

Next the seed 13(1) with b not zero, as explained in the previously.

$X_0$=(13+18)(1+18)= (31)(19)=**589** +(32+36(18))=1813 both iteration therefore Delta both=18(36+$a_{n-1}$)+ $N_{x-1}$,

32 is the a value from the first value 589. For the element after 1813=1813+(68+36)=(104)18+1813=3685$_1$

13(1+18)=**247**+18(13)=**481**+234=715+234=949

The subset generated for 13x1 is only iterated with 1+18, the left arm, if 1 remains stable produces the superset, as explained earlier, which is forbidden. Therefore we have no left arm, or 1,0, since we cannot allow multiplication by one to generate a non-prime subset.

The twins also have only two arms, 1,1, and 1,0, as 0,1=1,0;

7x7=**49** both=25x25=18(32)+49=**625**+(18(32+36))=1849 Delta=18($a_0$ +36)) + $N_{x-1}$

7+18=25(7)=**175**+18(7)=**301**+126=**427**+126=**553**+126=**679**+ 126= **805**

Total=2+4=6, with 175 taken out.

11x11=**121**+18(40)=841+18(40+36)=2209

11 steady=18(11)=**319**+198=**517**+198=**715**+198=913

From this list of calculated non-primes with the formulas discussed above, up to the first largest 1,1, iteration for LOPI 13 from the seed operations we calculated that we have 22 non-primes:

49, 85, 121, 175, 247, 265, 301, 319, 355, 391, 427, 445, 481, 517, 535, 553, 589, 625, 679, 697, 715, 805

from the seed operations, as shown above, which are easier calculated by binomial convolution for each seed, or (17 +18(a)(5+18(b)), 11+18(a)(11+18(b)), 7+18(a)(7+18(b)), with b≥a, (11+18(a)(11+18(b))) with b≥a. Leaving the other 23 elements generated by the linear progression in the LOPI 13 Superset, primes:

13, 31, 67, 103, 139, 157, 193, 211, 229, 283, 337, 373, 409, 463, 499, 571, 607, 643, 661, 733, 751, 769, 787 as the primes left in the LOPI 13 superset up to element $N_{45}$=805, 23 LOPI 13 Primes + 22 LOPI 13 NP=45. Theory with PNT, with pi(x)=N/ln N, gives roughly 120 primes and when we distribute these to the 6 known LOPI classes, it gives us about 20 primes per class. The real count is 23 for LOPI 13, which is very close to the theoretical value, even when excluding 2,3 from the theoretical count.

Given all these equations that sequentially and dependably populate all the non-prime subsets we can show the number of primes in each LOPI up to a chosen 18(a) value. The Delta that represents the changes in each binomial can be used to calculate the sequential non-primes. This method shows that in fact the primes exist where the factor generation mechanism can't produce factors that multiply to produce products, which is why primes are not the product of factors, since they are the gaps in each LOPI class. The distance between subsequent non-primes in each LOPI congruence class, (NP-NP/18)-1 equals the number of primes between subsequent non-primes, for example 125-35=90/18=5-1=4, or there are four primes between 125 and 35 in LOPI 17 residue class mod 18, which are 53, 71, 89, 107. Also the constant average of any LOPI class is a LOPI/LEPI value at any given element sum given by $a_x$ average for even number of elements, $LOPI_{n\ even}$, a odd, equals LEPI + 18($a_x$ average) where $a_x$ average = (LOPI element-LOPI/18 - 1) /2. For an odd number of elements, $LOPI_{n\ odd}$, a even, $a_x$ average equals LOPI + 18($a_x$ average), $a_x$ average = (a max LOPI element)/2. The average of the LOPI sets are alternating LOPI/LEPI values as the number of elements varies from odd/even, written by the linear binomials, 5+ 18(a) or 14 + 18(a), for example the sum up to LOPI 5 element 77, a = 4, therefore average LOPIx=4/2=2, therefore 5 +18(2)=41, 41 is the average value of all LOPI elements up to

element $X_{n=5}$ 77, or 5+23+41+59+77/5=41, or for even $X_{n=6}$ even for 95, (95-5/18 -1)/2 x average=4/2=2, therefore 14+ 18(2)=50,or  5+23+41+59+77+95/6=50.

The location of the non-primes is known by these equations and therefore the location and quantity of the primes,  LOPI Congurence Class – Non-prime calculated= Number of Primes in the integral wanted. The main issue here is to start the generation from the beginning and use the deltas in order as they build upon the seed operation and in which way the iteration of 18 occurred, Left, Right or Both.  This method can easily be programmed as well by a pretty simple algorithm that generates by the formulas and when the integer limit in question is reached, we count the elements batched in each sequence and then subtract from the total elements known in the LOPI superset.  Keep in mind that 2 and 3 are not included in the total when comparing to the theoretical values of pi(x) according to the prime number Theorem,  x/lnx, [24], there are 138 to be divided into 6 congruence LOPI classes = 23, exactly as found, at this point the primes are evenly distributed among the congruence classes.

N LOPI-LOPI/18 +1= Xn - Total elements in LOPI Superset – Total elements generated = Number of Primes in the LOPI superset.

The quantity calculated corresponds well to the theoretical value of pi(x)= N/ln N, which we then divide by 6 for the 6 congruence classes, since we assume even distribution.

Out of 45 elements LOPI 13 up to 805, we calculated 23 Primes and 22 Non primes,

Which (805/ln 805)/6= 20, which is close to the real value 23 primes, but still in the even distribution of primes and non-primes in the total 6 congruence classes.

The quantity of non-primes produced can be counted as we generate and then subtracted from the quantity known for each LOPI congruence class generated by the Dirichlet linear equations.  This will give us the exact number of primes and they are shown to always be present, as theory predicts for their even distribution in the congruence classes of the reduced residue system Mod 18.  This is a good way to generate the non-primes as well, but the factors are more elusive than by the convolution of the linear binomials and then matrix query.

The equations above can be used to sequentially generate the next non-primes in each subset, knowing that the primes are left out in the results of the convolution.

Every pair of possibilities in this LOPI P + LOPI P= LEPI template will be in the seed subset 13*7 and we can use our subsets of LOPI 1 to see if the product of these two LOPI elements is only once in subset (13 + 18(a)(7+18(b)), since this is a LOPI 1 non-prime generating subset, therefore we know it is a secondary prime and therefore meeting the Goldbach Conjecture. 13477*7=94339 is a prime pair, and therefore 13477+7 is a prime + prime from the LOPI templates of (13 + 18(a)+(7+18(b)), and it is in Row 1, Column 749, which corresponds to our earlier discussion that a=Row -1, and b=Column -1, or a=0, b=748, and the number of factors for the product 94339 is 2(x)+2, x being number of times in a subset of LOPI 1, also discussed in the above work, the +2 is for the trivial factors, 1 and itself, 2(1)+2=4, a secondary prime. In fact all the prime pair elements that "see-saw" to meet the prime + prime congruence relationship for LEPI 2, will be in only the one subset of their LOPI Superset 1, meaning that their product is a secondary prime and therefore is found only once in total in any of the subsets LOPI 1, with a prime from LOPI 7 and a prime from LOPI 13, which also could be checked as well from the original LOPI subsets of 7 and the LOPI subsets of 13, with no need to search the LOPI 1 subset, but one search of the specific non-prime LOPI 1, in this case, subsets for all secondary primes, elements occurring only one time in each of the 3 subsets, will give you Goldbach Prime + Prime solutions. All the templates for all the even congruence classes modulo 18 will be listed below.

We offer proof by construction but also by theory, for the Goldbach Conjecture to be false would require that, for example, the non-prime subsets of LOPI 1 that consists of the linear convolution of ((7 + 18(a))(13+18(b))or NP subset LOPI 1, consist only of tertiary primes or greater and therefore would require that the disjoint congruence classes of LOPI 7 and LOPI 13 not have evenly distributed primes or secondary primes from a=0 to b=c-1, that upon convolution produce secondary primes, as we know this is untrue, since clearly the secondary primes do exist as seen by the multiplicity of one in the data matrix tables, and as well the primes do exist as already by the Dirichlet Theorem for linear progressions using the totatives of mod 18. Also since LEPI 2=the Superset of LOPI 7 + the Superset of LOPI 13=LEPI 20, we must say that primes exist in these two linear

Dirichlet binomials that have been shown to yield convolution, secondary primes within the various subsets, the subsets are made of secondary, tertiary and higher primes, consecutively, we cannot make them without their predecessors, subsets cannot be made only of tertiary or greater non-primes, given that secondary primes will be infinite, as are the primes in each linear binomial in the convolution. The secondary primes occur at a pattern that is reflected in the prime element generation, as primes are made by the misses in the convolution of the linear LOPI equations with the variables a and b representing all the natural numbers from 0 to infinity. As well we know by the Dirichlet Theorem that the linear binomial equations with the LOPI, LOPI + 18(a) will contain an infinitude of primes as mod 18 is coprime to all LOPI, which equal the totatives of 18. And if primes are infinite then secondary primes are infinite and cannot be absent in the non-prime LOPI subsets. Furthermore, just like the point wise multiplication of the two linear binomials of LOPI 7, 13, results in a non-prime LOPI 1 product the point wise addition of LOPI 7 and LOPI 13, constitutes LEPI 2 + 18(a), or here since 2 is the LEPI we go to the next element, 2+18(1) = 20 =16+4 and 7+13.  All odd primes can be converted into their even LEPI counterpart and added to give their LEPI element mod 18.  Each LEPI congruence class mod 18 has more than one additive LOPI template of linear binomials from which to find the prime elements.  This is also without even using 2 or 3, which are not needed to solve the Goldbach Partition LOPI Seesaw as the gaps in the convolution show the primes already exist in time before an even LEPI element is generated via LEPI + 18(a). LOPI 7 + 18(a) + LOPI 13+ 18(b) = LEPI 2 + 18(c) , where a=b=c is a good first check. The LEPI elements between a=18(0) and b=(18(c-1) or (18(c)=b for LEPI >9, continue the linear progression at greater and greater lengths, producing primes and non-primes  as far as we wish to go, by the natural range of numbers a=0 to c, c-1. As there is an equal distribution of primes within these congruence classes, we can write each even number as a LOPI$_1$ Prime element + LOPI$_2$ Prime element, primes that come from the needed congruence classes of the LOPI reduced residue system of Mod 18.  To get the LEPI of the even number  there will always be a pair solution between a=0 and c, or c-1, that shows, as primes are evenly distributed gaps in the production of non-primes. Iteration by +18 and subsequent multiplication cannot produce all the natural odd numbers in the succession, but rather produces the non-primes only within each congruence class of each of the six LOPI. The gaps in the range of the domain of the convolution, the primes, add to create specific LEPI, according to the specific LOPI template that can be expressed as the congruence relations below, all LEPI have more than one Goldbach/Dirichlet LOPI class seesaw , thereby also providing all the secondary prime elements in the NP subsets.  The iteration and convolution of mod 18 is cyclic in the polynomial of the convolution, no prime can equal 0 in order to produce solutions, as the factor of 1 is not allowed.  Below are the specific LOPI linear Dirichlet binomials that create the Goldbach/Dirichlet LOPI See-saw  that yield all the even elements of the natural number line when the relationship of a to b  continues by a+1 and b -1 for each consecutive linear congruence. The elements in the even LEPI of modulo 18, written as linear equations, LEPI + 18(c)=(LOPI$_1$ +18(a) + (LOPI$_2$ + 18(b)).  Recently in conjunction with the Great Internet Mersenne Prime Search I received notification on April 14 that Oliveira and Silva have recorded that 3,325,581,707,333,960,528 is the smallest number that has no Goldbach partition with a prime. ( see Appendix).  A faster way to search possible solution is to use the Goldbach/Dirichlet LOPI seesaw algorithm to find the LOPI linear progressions that will work to yield solutions.  This integer is a LEPI 8 element, with c=1.8475453929633114e+17 and from template one, (13 + 18(a) (Mod 13 + 18(b), c/2≤b≤c, ((c-1)/2≤b≤c-1

 check multiplicity 1 in LOPI 7 sub 13*13

we use (8+ 18(1.8475453929633114e+17)≡13 + 18(0)mod 13+ 18(9.237726964816557e+16)  and by template two we see:

LEPI 8 +18(c)≡(1 + 18(a) (Mod 7 + 18(b) b≤c , a>0         Multiplicity 1check  LOPI 7     1*7 subset secondary

        LEPI 8 + 18(1.8475453929633114e+17)≡(1 + 18(1)(mod 7 + 18(1.8475453929633114e+17)

The first pair generated is (LOPI 1) 19 and  3,325,581,707,333,960,509 (LOPI 7)

and (LOPI 13) 13 and  3,325,581,707,333,960,515 (LOPI 13), obviously not a prime pair.  The algorithm can be coded to generate all the pairs meeting these two equivalence relationships and then the appropriate data matrix queried.

The range is large between $a_0=1$ and $b_{max}$ in which to find simultaneous primes in the linear progressions listed. For there to be no prime + prime elements in each linear progression is to say that there will be no prime element pairs in the linear Dirichlet progressions between 1 + 18(1) to (1 + 18( 1.8475453929633114e+17) and 7 + 18( 1.8475453929633114e+17) to 7 + 18(0)., nor a pair of simultaneous primes in the other template, 13*13 as well. The specific Goldbach/Dirichlet LOPI Seesaw progressions can be calculated and then the elements individually tested for primality by zero presence in the subsets or the products testing for multiplicity 1 in the known subsets, 13*13 subset LOPI superset 7 and subset 7*1 in LOPI superset 7. This method will more efficiently find and confirm prime elements that will solve the Goldbach Conjecture for the integer 3,325,581,707,333,960,528. When looked at as the range of the linear equations that need to just have one pair prime at the same time between a miniumu and b max, it looks even more unlikely that these primes do not exist, as we know they are distributed throughout the whole LOPI congruence class to infinity as they are reduced residue Dirichlet linear binomial progressions mod 18.

Solutions are the secondary prime in the LOPI 11 subset 13*5, or elements appearing only once, or a search for elements that appear 0 times in any subset for the respective LOPI superset 5 and LOPI superset 13 can also yield the Goldbach/Dirichlet pair of elements.

Goldbach/Dirichlet LOPI seesaw templates for all even natural numbers ≥6 residue mod 18, not to include 3+3=6. LEPI residues mod 18=0, 2, 4, 6, 8, 10, 12, 14, 16. For all values using LOPI 1, a>0.

LEPI 0 +18(c)≡(1 + 18(a) (Mod 17 + 18(b),  $a_0=1$, $b_{max}=c-1$ LEPI>9,   LOPI 17 1*17 subset sec.

    Continue the see-saw for each pair LOPI + LOPI, as (a) increases by 1, (b) decreases by 1,

(LOPI + LOPI) ≡ LOPI + 18($_{a0,a+1..a=bmax,}$)(mod LOPI +18($_{b\ max,b-1..b=a0}$)), b max defined in each as listed below.

LEPI 0 +18(c)≡(1 + 18(a) (Mod 17 + 18(b),  $a_0=1$, b≤c-1     LOPI 17 1*17 subset sec.

        ≡(7 + 18(a) (Mod 11 + 18(b), $a_0=0$, b≤c-1     LOPI 5 7*11 subset secondary

        ≡(5 + 18(a) (Mod 13 + 18(b), $a_0=0$, b≤c-1     LOPI 11 13*5 subset secondary

LEPI 2 +18(c)≡(1 + 18(a) (Mod 1 + 18(b), $a_0=1$, c/2≤b≤c even, b=(c-1)/2≤b≤c odd,   LOPI 1 (1*1) sub

$\equiv (7 + 18(a)$ (Mod $13 + 18(b)$ $a_0=0$, $b \leq c-1$  LOPI 1 13*7 subset sec.

LEPI 4 $+18(c) \equiv (5 + 18(a)$ (Mod $17 + 18(b)$, $a_0=0$, $b \leq c-1$,  LOPI 13 17*5 sub sec

$\equiv (11 + 18(a)$ (Mod $11 + 18(b)$, $a_0=0$, $b=c/2 \leq b \leq c$ even, $(c-1)/2 \leq b \leq c-1$, odd  LOPI 13 11*11 sub

LEPI 6 $+18(c) \equiv (1 + 18(a)$ )(Mod $5 + 18(b)$, $a_0=1$, $b \leq c$,  LOPI 5 1*5 subset secondary

$\equiv (7 + 18(a)$ (Mod $17 + 18(b)$, $a_0=0$, $b \leq c-1$  LOPI 11 7*17 subset secondary

$\equiv (11 + 18(a)$ (Mod $13 + 18(b)$, $a_0=0$, $b \leq c-1$  LOPI 17 13*11 subset secondary

LEPI 8 $+18(c) \equiv (1 + 18(a)$ (Mod $7 + 18(b)$ $a_0=1$, $b \leq c$  LOPI 7 1*7 subset secondary

$\equiv (13 + 18(a)$ (Mod $13 + 18(b)$, $a_0=0$, $c/2 \leq b \leq c$, $((c-1)/2 \leq b \leq c-1$, odd) LOPI 7 sub 13*13

LEPI 10 $+18(c) \equiv (5 + 18(a)$ (Mod $5 + 18(b)$, $a_0=0$ $b \leq c/2$, even, $c-1 \leq b \leq c$ odd LOPI 7 5*5 subset

$\equiv (11 + 18(a)$ (Mod $17 + 18(b)$ $a_0=0$, $b \leq c-1$  LOPI 7 11*17 subset sec.

LEPI 12 $+18(c) \equiv (1 + 18(a)$ (Mod $11 + 18(b)$, $a_0=1$, $b \leq c$  LOPI 11 1*11 subset secondary

$\equiv (5 + 18(a)$ (Mod $7 + 18(b)$, $a_0=0$, $b \leq c$,  LOPI 17 5*7 subset sec

$\equiv (13 + 18(a)$ (Mod $17 + 18(b)$ $a_0=0$, $b \leq c-1$  LOPI 5 13*17 subset secondary

LEPI 14 $+18(c) \equiv (1 + 18(a)$ (Mod $13 + 18(b)$, $a_0=1$, $b \leq c$,  LOPI 13 1*13 subset

$\equiv (7 + 18(a)$ (Mod $7 + 18(b)$, $c/2 \leq b \leq c$, $(c-1)/2 \leq b \leq c$, odd  LOPI 13 7*7 subset

LEPI 16 $+18(c) \equiv (5 + 18(a)$ (Mod $11 + 18(b)$, $a_0=0$, $b \leq c$,  LOPI 1 5*11 subset sec.

$\equiv (17 + 18(a)$ (Mod $17 + 18(b)$ $c/2 \leq b \leq c-1$, $(c-1)/2 \leq b \leq c-1$, odd  LOPI 1 17*17 subset

LEPI 88 $\equiv (5 + 18(a)$ (Mod $11 + 18(b)$ $a_{min}=0$, $b \leq c$, $c=4$

| $a_i=0$ | 5 | 83 | $b_{i=c}$ 4 | $I_a$ a0, 5 (a+1) ↓ |  | 11 | $F_b$ b=a0 |
|---|---|---|---|---|---|---|---|
| 23 |  | 65 |  | 23 +1 ↓ |  | 29 | ↑ -1 |
| 41 |  | 47 |  | 41 +1 ↓ | See-Saw | 47 | ↑ -1 |
| 59 |  | 29 |  | 59 +1 ↓ |  | 65 | ↑ -1 |
| $a_f$ max, 77 | 11 | $b_f=0$ |  | $F_a$ a=b max 77 |  | $I_b$ $b_i$ max 83 ↑b-1 | |

The coefficient *a initial* always equals 0 and *a final* will be equal to *b initial*. I have listed all the cases to be clear for the research example. In reality the rules can be simplified for the algorithm based on the even integer one needs to inquire about as different additive cases vary a small amount, such as when LOPI 1=LOPI 2 the relationship of *b to c* changes. Although this list is long the actual equation is simple and to get b is not difficult nor to understand if it is equal to c or not. All even integers can be represented by many odd additive pairs from the Dirichlet linear progressions that make the reduced residue system Mod(18). As it is a reduced residue system, as discussed above, there are infinite cyclic primes distributed within these linear progressions. The Goldbach Partition solutions, using the infinite congruence classes into which all prime odd integers exist, LOPI residue mod 18, are possible to calculate and more than one solution will exist. As Christian Goldbach conjectured, we find it possible and offer this as numeric proof that a system does exist that explains why his conjecture is true. The Goldbach Partition solutions can be generated using the Dirichlet linear equations along with the Euler function totatives, $\varphi(18)$, as the residues coprime to mod 18, and making templates based on how the LOPI residues combine to form the complete set congruence classes of the even residues of modulo 18. The LOPI/LEPI relationship allows us a core integer characteristic by which to understand how to use the LOPI reduced residue classes in order to make a predictable see-saw partition congruence relation. It is a simple additive operation. We don't have to go through all the integers below 28, just the LOPI residue of modulo 18 as they are already prime and are the coprime totatives to mod 18. The linear Dirichlet Progressions will afford all possible LOPI prime + LOPI Prime that yield LEPI 28. 28 is a LEPI residue 10 element, 5+5, but also we must use 8+2 or 17+11, as to make use of all the LOPI >18 we need to take advantage of all the LOPI that add together to give the LEPI, therefore 17 and 11 are LEPI/LOPI equivalents of 8 and 2. Template 17,11 pair yields: $(10 + 18(1) \equiv 11 + 18(0)_{(mod\ 17\ +\ 18(0))}$ and 5,5 pair, yields 23 +5 = 28, as $(10 + 18(1) \equiv 5 + 18(0)_{(mod\ 5\ +\ 18(1))}$. The relations can then be manipulated directly to create the whole set of even integers if we wish by going consecutively by +1 in the residue linear equation and by -1 in the modulo linear equation, as decided by the beginning value of b based on c, as shown above to create the full set of pairs for the Goldbach/Dirichlet LOPI congruence seesaw.

LOPI mod 18, based on the mechanism of iteration and factor generation through convolution can partition the congruence classes that can be paired as shown above to solve the Goldbach/Dirichlet congruence LOPI seesaw. We can query the subsets of the product LOPI for multiplicity in the integral listed in the algorithm section to see if it is secondary or to each template for secondary primes or of the individual elements, as when the pair is multiplied, if both are prime the product will be found in the specific non-prime subset predictable by multiplying the two LOPI in the equation. For instance above we had 23 + 5=28. If these secondary primes did not occur cyclically then we would have trouble filling in the subset elements that are > secondary prime, as each degree of primality originates from the ones before it, as 0-Prime, 1, 2, 3. The elements of the primes, 1, 2, 3... are created through the sequential iteration of 18 as discussed above and the primes are skipped regularly in the linear convolution of LOPI modulus 18. Due also to the basic mathematical predictions that are possible due to the constant core LOPI identity being equal to the totatives of mod 18, the grouping identity withstands these operations and proves to be a very predictable and important tool for partitioning and researching the distribution of integers by their degree of primality. Using this LOPI reduced residue mod 18 system we can predict and create many different patterns, as needed. The Goldbach Partitions in the form of $P_{LOPIa + 18(a)} + P_{LOPIb + 18(b)} = N_{LEPIc + 18(c)}$, can be generated for all even integers, LEPI mod 18 congruence classes, within this reduced residue system of mod 18, which contains all the odd natural numbers capable of being prime, we have an infinite source of primes and non-primes being generated only through the Dirichlet linear progressions, while the non-prime elements are equal to the products generated through the specific convolution of the same binomial equations making the complete set of non-prime elements within each LOPI, leaving the complement of these non-prime subsets equal to the prime subset of each LOPI. The empty set is the total LOPI subsets – Prime subset = 0. The relationship LEPI + 18(c)≡LOPI$_1$ + 18(a) (mod LOPI$_2$ + 18(b)) finds all possible prime, prime solutions in the integral of $a_0$ and $b_{max}$ and as they originate from different distances from zero, their prime/non-prime pattern will never be able to be only NP + NP, NP + P, P+ P or P+NP, they must contain the four states these pairs can exist in and can't stay in any one pair state infinitely or we would stop using each LOPI residue and as the linear equation clearly shows counting by 18 from the origin of the residues completely and sequentially creates the complete set of elements in each residue class. None can be missed in the modular math of modulus iteration. Also the probability is zero of having a prime in the subsets since we allow NP products only in the form of factors other than one, and since by definition primes are only formed in multiplication by 1, there can be no primes in the NP subsets. Since we have sequentially created each element in each congruence class using the linear Dirichlet Binomial with mod 18, and then use these same linear binomials to create the factor pair/product, we know primes are the gaps in the convolution, as the equivalence relation discussed above maps over only to the non-primes with factors >1 since 1 is not allowed as a factor, in the NP polynomial generators; Seed + 18(LOPI)(b) + 18(LOPI)(a) + 324(ab) can only subtract from it's non-prime linear value in the LOPI superset to equal zero. It follows logically also that the seed constant in these polynomials are formed only of the linear Dirichlet progressions, (LOPI + 18(a)(LOPI + 18(b)) with two LOPI residues that are not 1, except for the case of the identity NP subset that includes 1 as its LOPI, but in these cases we never allow b to equal zero, therefore 1 is always added to a quantity greater than zero to make the factor in order to avoid any prime generation.

The Goldbach Partition solutions using the infinite congruence classes into which all prime odd integers exist will provide sufficient primes to produce their corresponding even residue classes. The iterative factor generation mechanism of the LOPI reduced residue system modulo 18 uses linear Dirichlet progressions to generate factors and then convolution of these vector and row columns to produce all the non primes coprime to 18 partitioned into well defined subsets by realizing the constant core lowest odd partition identity can be used to generate and group the non-primes, leaving us with proper NP subsets that upon subtraction from the larger congruence class superset, leaves only the prime subset:

{Superset A}\ NP {sub B, sub C, sub D}={NP'}={Prime subset} given the probability is zero that a prime can be in the subsets, {NP-NP'}=Ө or {NP-P}=0 for each LOPI reduced residue class. The attached convolution algorithm queries multiplicity in a specific integral in the data matrix subsets as a Deterministic Primality test, where 0 not found indicates a Prime and found at one location indicates a secondary prime. A secondary that has the product of LOPI residue 5 is important to consider a secondary since it is the only non-prime ending in a 5 that will count to a multiplicty of one. These products are important to be noted as secondary, since in the Goldbach Conjecture 5 is one of the prime elements in the adding of the elements of the specific LOPI Dirichlet linear equations.  Any other product with a five in it will not have a multiplicity of 1, as the factor will always have other factors since it is not the LOPI residue itself, 5. We then add in another option to obtain the factors of the integer using the found feature, row and column, by which we know also the factors, as row and column are b+1 and a+1 in the case of all but LOPI 1,1 subset as discussed above, where m=a and n=b since no multiplication by one is ever allowed. Please note the LOPI$^3$ unique example of multiplicity of one as listed to be excluded up front as explained earlier in this work.

**Conclusion**

The degree of primality of an integer can be determined by its multiplicity in a finite set of sequential infinite non-prime element data matrices disjoint from the other LOPI non-prime subsets, uniquely grouped and generated under their LOPI family constant identity using the LOPI reduced residue system of mod 18 and the restrictions built in as discussed above. By viewing each natural integer as an unrestricted partition of the lowest odd residues mod 18 we are able to group the odd natural integers into disjoint super-sets upon this lowest odd partition identity, contained in all integers. The core of the elements in these classes dictates what will occur to the integers in the classes upon certain operations, such as shown for the Sophie Prime Germain operation, 2p + 1.  The effect of the operation carries through to all the elements within the classes as they are the LOPI core with sequential degrees of the neutral cycle of 18 iterated to infinity.  The neutral value in the Dirichlet linear progressions advance the LOPI to the next partition while maintaining the core LOPI identity. The same holds true for the non-prime generating elements through the convolution of specific linear binomials, as shown above. The unique factors and subsequent product elements generated through consecutive iteration and convolution populate the proper non-prime subsets of each LOPI mod 18 reduced residue class, yielding the complete 6 LOPI disjoint reduced residue classes mod 18, each superset {A} is equal to {NP Sub B} + {NP Sub C} + {NP Sub D} + {Prime Subset}, remembering that each LOPI Superset has 3-4 non-prime subsets, here we just list three, B,C,D, for the three subset LOPI congruence classes mod 18.

With the LOPI system of binomial convolutions, we are generating the non-primes and their factors, rather than inquiring as to whether a number is divisible by certain numbers or how many factors. Primality multiplicity = 0 for the return result of number of non-prime LOPI subsets in which we find the LOPI integer, 0, prime, 1 is a secondary prime, 2 a tertiary and so forth as discussed with the relationship, 2 + 2(n), with n=numbers of times number(n) is found in LOPI N subset, the result equals the number of factors of n, with primes never being found in a subset, therefore n=0 and primality=0, and number of factors = 2, trivial. Using the attached algorithm we believe the LOPI mod 18 system of partitioning and subsequent query to be a polynomial time efficient method for deterministic primality determination and factor generation as well.

We can work in the negative integers as well by crossing over left to zero and the functions change in that the negative lowest odd partition becomes the LOPI -18. For example, LOPI 5-18=LOPI -13, and vice versa, [13,-5] therefore LOPI 5 + 18(a) =-13 LOPI + 18(a+1). (LOPI$_1$ + 18(a))+(-LOPI $_2$+ 18(a+1)=0, 1/23+23/1=0, (5 + 18(a)+(-13 + 18(a+1)=0. We can approach a space in which an integer exists from a negative LOPI inverse and/or the positive LOPI, for example to arrive at LOPI 5 element, 617 in the linear Dirichlet binomial progressions, we can find it by 5 + 18(34)=617 and we can also generate it by using -13 LOPI by using (-13 + 18(a+1)=617. The same holds true for [7, -11] and [1,-17] the reverse. If we want to use the negative system of the mod 18 then we negate both elements of the linear binomial equation, -LOPI + (-18(a), for example (-13 + (-18(a)), and upon convolution we yield the same positive non-prime LOPI element expected. (-13 + -18(a)(-7 + (-18(b)) will yield a positive LOPI 1 non-prime element at the same cycle of 18 a and b. Using a=1 as an example, we get (-13 + -18(a)(-7 + (-18(b)) = − 31( 25)=775, which equals (13 + 18(1))(7 + 18(1)). And as we do not allow the multiplication by 1,-1, or (LOPI 1 +18(0), (not allowed), the primes are also left in the negative system of mod 18 upon using the convolution method as their inverses are not created, as discussed above, by the LOPI non-prime convolution of their specifically needed LOPI linear equations. The only way they could be put back into the positive number system would be by multiplication by -1, which is the main requirement of this system, no convolution by 1 in the Linear binomial factor generator of (LOPI 1 + 18(a)), when it is part of the convolution, the coefficient to 18 with LOPI 1 can never be allowed to be 0, only >0.  This is one of the main reason the system works to partition the prime from the non-prime elements with non-trivial factors, factors other than 1, the only factor of a prime besides itself.   The fact that all the non-primes coprime to 18 are generated from the first few seed factors coprime to 18 in each LOPI reduced residue class mod 18, is elegant and deterministic at the same time.  Instead of 1 pattern containing the primes, we see 6 patterns exist as determined by the LOPI core identity mod 18, where LOPI/LEPI  N($_{a,bc..n}$)= (a + b + c...n).,= N LOPI/LEPI mod 18.  As Euclid called the primes, linear elements, we see these core linear LOPI progressions, advance by 18(a), 18(b) iteration and convolution to create the plane and solid elements.  The reduced residue system modulo 18 combined with linear convolution, ELMA-18,  provides a possible solution to the P NP computational problem, providing a polynomial ask and a deterministic polynomial answer to the question: Is an integer a prime/degree of primality or a non-prime and if non-prime, what are the non-trivial factors. The mathematical model shows that the factors are all generated in polynomial time upon the convolution of the specific linear Dirichlet progressions and upon matrix subset query within the already known congruence family class based on the LOPI of the integer, factors of any large integer are designed to be the index vector and row quantities and therefore reported immediately upon the location of the integer in the subset matrix.  A Zero find reveals a prime, they cannot map into the non-prime subsets as a row/vector factor of one is never allowed in the convolutions. Primes are the only numbers to have only trivial factors, therefore an excellent partitioning parameter for prime and non-prime numbers. All prime elements are partitioned out as they are the compliment to the non-prime subset elements within each superset congruence class.

We used Excel 2013 to generate the LOPI supersets by (LOPI + 18(a)), a≥0 and to model the LOPI mod 18 Dirichlet linear binomial convolution system, to generate the unique non-trivial factors and subsequent non-prime products for each LOPI non-prime subset of the 6 LOPI congruence classes mod 18. The Excel model is available upon request. laurellynn.mcclure@studio.unibo.it

**Acknowledgments**

I would like to thank my family for all the time and support they lovingly gave me and dedicate this work to them: Maurizio, Maxwell, Lucia and Emma.
As well to Dr. John Stickney, Distinguished Research Professor, Department of Chemistry, University of Georgia, Athens, a constant and steady mentor through all the years I have known him.
and to my good friend, and to Dr. Allen J. Bard, UT Austin, the real deal.
Dr. Constance Mwangemi, who sings too, since we know, like Sophie Germain and all the others, that….
" A bird doesn't sing because it has an answer, it sings because it has a song." Joan Walsh Anglund
Thank you to the state of NC and the excellent Public School and University system and to my trusty Hurbert , my Texas Instrument TI-36X Pro.

**REFERENCES**

[1] Beutelspacher, Albrecht, *Numbers*, Mineola, New York, Dover, (2015), pg. 11-14.    Originally published as C.H. Beck oHG, *Zahlen: Geschichete, Gesetze, Geheimnisse*,   Munich, Germany (2013).

[2] Zaiger, D., *Newman's Short Proof of the Prime Number Theorem*, The American Mathematical Monthyly, Vol. 104, No. 8 (Oct., 1997), pp. 705-708.

[3] Sabia, J., & Tesauri, S. (2009). The Least Prime in Certain Arithmetic Progressions. *The American Mathematical Monthly,116*(7), 641-643. Retrieved from http://www.jstor.org/stable/40391172

[4] Ramanujan, Sirindudi, *A proof of Bertrand's Postulate*, Journal of the Indian Mathematical Society, XI, 1919, pp. 181-182.

[5]  Euclid, *The Elements,* Book 7 Vol. 2, Dover Edition, Cambridge University Press, (2017),    pg. 284-285
The Thirteen Books of Euclid's *Elements*, Infinity of Primes Vol.2 The Elements,   Translated by Sir Thomas L. Heath, Honorary Fellow of Trinity College, Cambridge, Dover Publication (Three Volumes), New York 1956

[6] George E. Andrews,  *Number Theory*, 2nd edition, New York, Dover Publications (1994), pp. 54-55.

[7] Erdos, P., *On a problem concerning systems of congruences*, Mat. Lapok 3(1952), 122-128.

[8] Waid, Carter, *On Dirichlet's Theorem and Infinite Primes*, Proceedings of the American Mathematical Society, Vol. 44, No. 1, (May 1974), pp. 9-11.

[9] H. Maier, *Primes in short intervals*, Michigan Math. J. 32 (1985), pp. 221–225.

[10]  SHIU, D., & SHIU, P. (2011). A lower bound for the prime counting function. *The Mathematical Gazette, 95*(534), pp. 433-436. Retrieved from http://www.jstor.org/stable/23248511

[11] George E. Andrews,  *Number Theory*, 2nd edition, New York, Dover Publications, (1994), pg. 13.

[12] Allenby, R.B.J.T., *Rings, Fields and* Groups, 2nd edition (2001), Oxford, UK and Burlington MA, Butterworth-Heinemann Publishers, (1991 first published) pp. 202-220.

[13] Pinter, Charles C., *A Book of Abstract Algebra,* 2nd edition, New York, Dover Publications, (2016), © 1982, 1990 by Charles Pinter. pg. 25-32.

[14] Allenby, R.B.J.T., *Rings, Fields and* Groups, 2nd edition (2001), Oxford, UK and Burlington MA, Butterworth-Heinemann Publishers, (1991 first published) pp. 60-61.

[15] Barnabei, Marilena and Bonetti, Flavio, *Elementi di Aritmetica Modulare*, 2nd edition (2014), Bologna, Italy, Societa' Editrice Esculapio, pp. 21-26.

[16] George E. Andrews, *Number Theory*, 2nd edition, New York, Dover Publications, (1994), pg. 60-61.

[17] Clapham, C., Nicholson, J., *Oxford Concise Dictionary of Mathematics*, Oxford, UK, Oxford Univ. Press, (2014), pg. 115.

[18] Amy Dahan Dalmédico, *Sophie Germain,* Scientific American, Vol. 265, No. 6 (DECEMBER 1991), pp. 116-123.

[19] George E. Andrews, *Number Theory*, 2nd edition, New York, Dover Publications, (1994), pp. 113-116.

[20] Indlekofer Karl-Heinz, Ja'rai, Antal, *Largest Known TwinPrimes and Sophie Germain Primes*, Mathematics of Computation, Volume 68, Number 227, pg. 1319, part 4.

[21] Wanko, Jeffrey J., *The Legacy of Marin Mersenne: The Search for Primal Order and the Mentoring of Young Minds Author(s),* The Mathematics Teacher, Vol. 98, No. 8 (APRIL 2005), pp. 525-529 Published by: National Council of Teachers of Mathematics Stable URL: https://www.jstor.org/stable/27971800 Accessed: 17-05-2019 06:57 UTC

[22] Jenkins, Paul, *Odd Perfect Numbers Have a Prime Factor Exceeding 107*. *Mathematics of Computation,* (2003), *72* (243), 1549-1554. Retrieved from http://www.jstor.org/stable/4099849

[23] Putnam, T. (1910), *Perfect Numbers*, The American Mathematical Monthly*, 17*(8/9), 165-168. doi:10.2307/2973855

[24] Sandifer, C. Edward, *The Early Mathematics of Leonard Euler* , The Mathematical Association of American, 2007, pg. 253.

[25] ArXiv:0708.2564v3 [math.HO], Bell, Jordon, 2017, pg.2.

[26] Goles, Eric, Schulz, Oliver and Markus, Mario, *A Biological Generator of Prime Numbers*, Nonlinear Phenomena in Complex Systems, 3:2 (2000) pp.208-213.

[27] Ferdeghini, Andrea, Franchini, Anna, (editor), *I Numeri Primi,* Milano, Italy, RBA, Mondo Mathematico, (2018), pp. 58-59.

Figure 1 Number line reflecting grouping identity constant, linear LOPI + 18(a), or $N_{LOPI} \equiv LOPI$ (Mod 18)

LOPI = φ(18); 1, 5, 7, 11, 13, 17, to reveal seed LOPI factors for non-prime subset generation in Fig.2

[LOPI=a+b+c..= $N_{LOPI} \equiv LOPI$ (Mod 18)]

| 18 *,+,neutral absorbing additive | LOPI 1 | LOPI 11 | LOPI 13 | LOPI 5 | LOPI 7 | LOPI 17 |
|---|---|---|---|---|---|---|
|  | **1** | **11** | **13** | **5** | **7** | **17** |
| 18 | 19 |  |  | 23 | **25** |  |
|  |  | 29 | 31 |  |  | **35** |
| 36 | 37 |  |  | 41 | 43 |  |
|  |  | 47 | **49** |  |  | 53 |
| 54 | **55** |  |  | 59 | 61 |  |
|  |  | **65** | 67 |  |  | 71 |
| 72 | 73 |  |  | **77** | 79 |  |
|  |  | 83 | **85** |  |  | 89 |
| 90 | **91** |  |  | 95 | 97 |  |
|  |  | 101 | 103 |  |  | 107 |
| 108 | 109 |  |  | 113 | 115 |  |
|  |  | **119** | **121** |  |  | 125 |
| 126 | 127 |  |  | 131 | 133 |  |
|  |  | 137 | 139 |  |  | **143** |
| 144 | 145 |  |  | 149 | 151 |  |
|  |  | 155 | 157 |  |  | 161 |
| 162 | 163 |  |  | 167 | **169** |  |
|  |  | 173 | 175 |  |  | 179 |
| 180 | 181 |  |  | 185 | **187** |  |
|  |  | 191 | 193 |  |  | 197 |
| 198 | 199 |  |  | 203 | 205 |  |
|  |  | 209 | 211 |  |  | 215 |
| 216 | 217 |  |  | **221** | 223 |  |
|  |  | 227 | 229 |  |  | 233 |
| 234 | 235 |  |  | 239 | 241 |  |
| 243 |  | 245 | 247 |  |  | 251 |
|  | 253 |  |  | 257 | 259 |  |
| 261 |  | 263 | 265 |  |  | 269 |
|  | 271 |  |  | 275 | 277 |  |
| 279 |  | 281 | 283 |  |  | 287 |
|  | **289** |  |  | 293 | 295 |  |
| 297 |  | 299 | 301 |  |  | 305 |

Figure 2 – Universal set = Total of Super-sets {LOPI 1} + {LOPI 5}…{LOPI 17}, LOPI 1,5,7,11,13,17
  LOPI Superset = Non-prime subsets,({$B_{lopi}, C_{lopi}, D_{lopi}$…}) + {$Prime_{LOPI}$ subset}
   Proper Non-Prime Subsets seeds in blue for convolution of the LOPI Dirichlet linear equations
   Ex. LOPI 17 subsets=(LOPI 17 + 18(a))(LOPI 1+18(b), a≥0, b>0
   ((LOPI 13 + 18(a))(LOPI 11+18(b) and ((LOPI 7 +18(a))(LOPI 5+18(b), a≥0, b ≥0
   ($LOPI_1$ + 18(a))($LOPI_2$ + 18(b)) with a≤b≥0, when LOPI $1^2$ sub LOPI 1, a≤b≥1

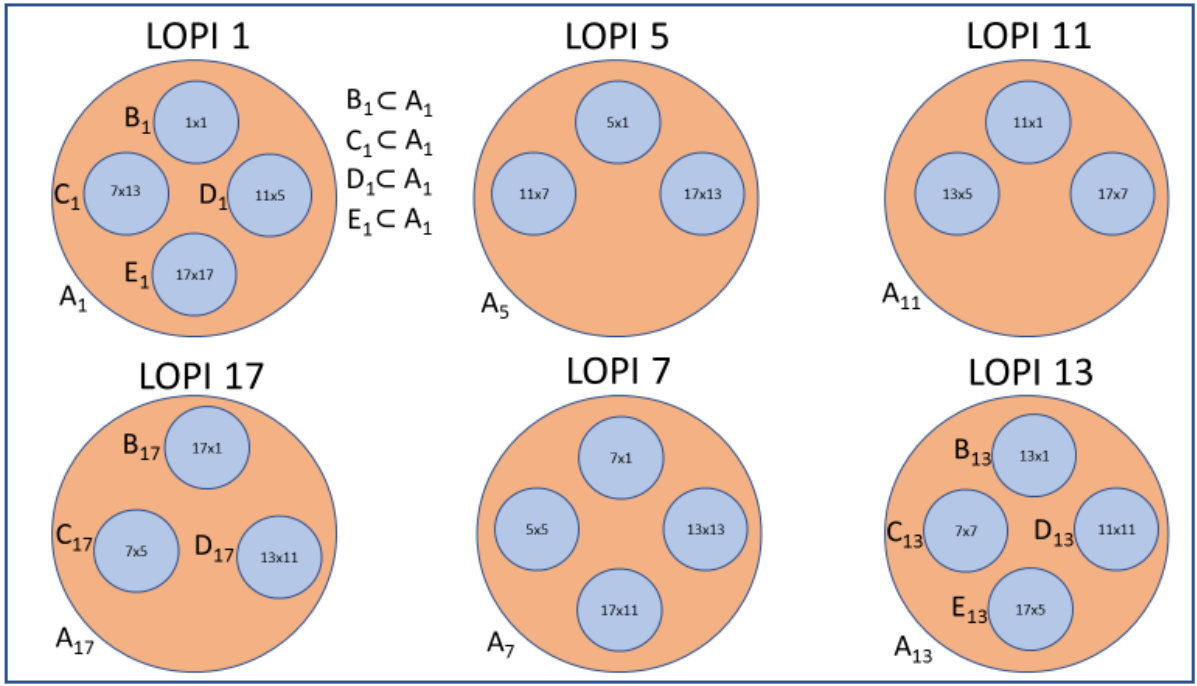

6 LOPI: 1, 5, 7, 11, 13, 17

Figure 3

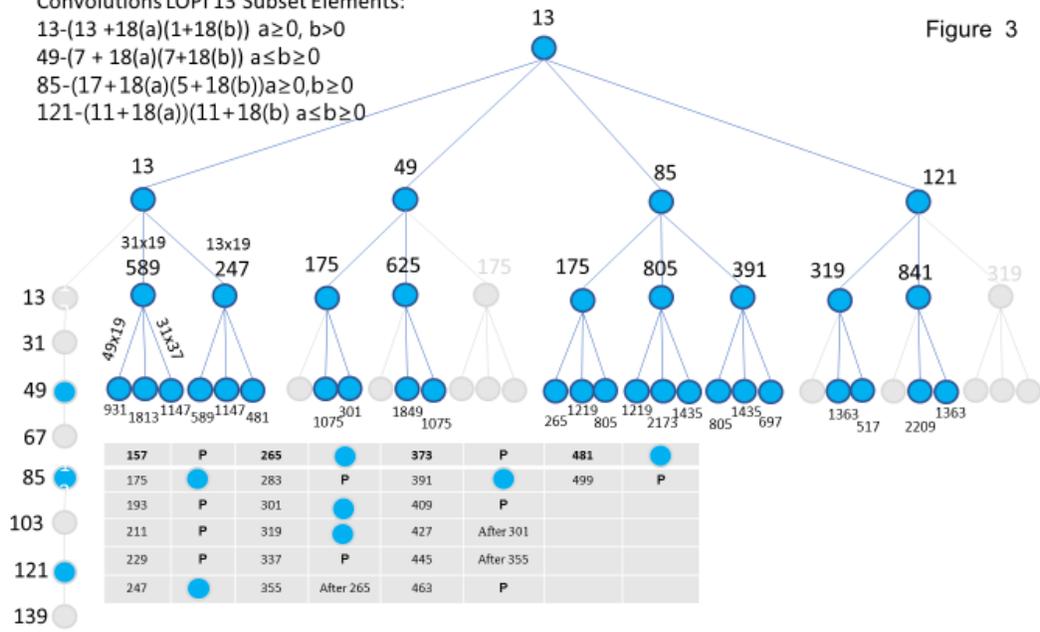

Figure 4 –
Convolution example of LOPI NP proper subset $B_{17}$ – $((LOPI\ 17+18(a))((1+18(b)),\ a\geq0,\ b>0$
For Superset LOPI 17, $A_{17}$, generated by the Dirichlet Linear Progression- LOPI 17 + 18(a)

|    | 1       | 2    | 3     | 4     | 5     | 6     | 7     | 8     |
|----|---------|------|-------|-------|-------|-------|-------|-------|
| 0  | Factors | 19   | 37    | 55    | 73    | 91    | 109   | 127   |
| 1  | 17      | 323  | 629   | 935   | 1241  | 1547  | 1853  | 2159  |
| 2  | 35      | 665  | 1295  | 1925  | 2555  | 3185  | 3815  | 4445  |
| 3  | 53      | 1007 | 1961  | 2915  | 3869  | 4823  | 5777  | 6731  |
| 4  | 71      | 1349 | 2627  | 3905  | 5183  | 6461  | 7739  | 9017  |
| 5  | 89      | 1691 | 3293  | 4895  | 6497  | 8099  | 9701  | 11303 |
| 6  | 107     | 2033 | 3959  | 5885  | 7811  | 9737  | 11663 | 13589 |
| 7  | 125     | 2375 | 4625  | 6875  | 9125  | 11375 | 13625 | 15875 |
| 8  | 143     | 2717 | 5291  | 7865  | 10439 | 13013 | 15587 | 18161 |
| 9  | 161     | 3059 | 5957  | 8855  | 11753 | 14651 | 17549 | 20447 |
| 10 | 179     | 3401 | 6623  | 9845  | 13067 | 16289 | 19511 | 22733 |
| 11 | 197     | 3743 | 7289  | 10835 | 14381 | 17927 | 21473 | 25019 |
| 12 | 215     | 4085 | 7955  | 11825 | 15695 | 19565 | 23435 | 27305 |
| 13 | 233     | 4427 | 8621  | 12815 | 17009 | 21203 | 25397 | 29591 |
| 14 | 251     | 4769 | 9287  | 13805 | 18323 | 22841 | 27359 | 31877 |
| 15 | 269     | 5111 | 9953  | 14795 | 19637 | 24479 | 29321 | 34163 |
| 16 | 287     | 5453 | 10619 | 15785 | 20951 | 26117 | 31283 | 36449 |

Below, portion of the data matrix of the convolution for non-prime subset 17*1 of Superset LOPI 17

|      |         | 556        | 557        | 558        | 559        | 560        | 561        |
|------|---------|------------|------------|------------|------------|------------|------------|
|      | Factors | 9991       | 10009      | 10027      | 10045      | 10063      | 10081      |
| 8822 | 158795  | 1586520845 | 1589379155 | 1592237465 | 1595095775 | 1597954085 | 1600812395 |
| 8823 | 158813  | 1586700683 | 1589559317 | 1592417951 | 1595276585 | 1598135219 | 1600993853 |
| 8824 | 158831  | 1586880521 | 1589739479 | 1592598437 | 1595457395 | 1598316353 | 1601175311 |
| 8825 | 158849  | 1587060359 | 1589919641 | 1592778923 | 1595638205 | 1598497487 | 1601356769 |
| 8826 | 158867  | 1587240197 | 1590099803 | 1592959409 | 1595819015 | 1598678621 | 1601538227 |
| 8827 | 158885  | 1587420035 | 1590279965 | 1593139895 | 1595999825 | 1598859755 | 1601719685 |
| 8828 | 158903  | 1587599873 | 1590460127 | 1593320381 | 1596180635 | 1599040889 | 1601901143 |
| 8829 | 158921  | 1587779711 | 1590640289 | 1593500867 | 1596361445 | 1599222023 | 1602082601 |
| 8830 | 158939  | 1587959549 | 1590820451 | 1593681353 | 1596542255 | 1599403157 | 1602264059 |
| 8831 | 158957  | 1588139387 | 1591000613 | 1593861839 | 1596723065 | 1599584291 | 1602445517 |
| 8832 | 158975  | 1588319225 | 1591180775 | 1594042325 | 1596903875 | 1599765425 | 1602626975 |
| 8833 | 158993  | 1588499063 | 1591360937 | 1594222811 | 1597084685 | 1599946559 | 1602808433 |
| 8834 | 159011  | 1588678901 | 1591541099 | 1594403297 | 1597265495 | 1600127693 | 1602989891 |
| 8835 | 159029  | 1588858739 | 1591721261 | 1594583783 | 1597446305 | 1600308827 | 1603171349 |
| 8836 | 159047  | 1589038577 | 1591901423 | 1594764269 | 1597627115 | 1600489961 | 1603352807 |

**Appendix**: New Algorithm for efficient deterministic Primality testing and factor determination

ELMA-18:

Possibly in P time [15, pg. 107-110]  Prime=0=yes , >0=no list factors

**1**-(Sum digits excluding nine or divide by 18 and get LOPI residue)

**2**-Convolute appropriate subsets by linear convolution of $(LOPI_1+18(a))((LOPI_2 + 18(b))$

$a \geq$, $b \geq 0$  appears to O(x), linear, f(a,b) both linear time,

when $LOPI_1=LOPI_2$, ≈1/2 size of whole subset, $a \leq b \geq 0$ poly becomes:

$((LOPI_1 + 18(b-1)(LOPI_1+18(b))) =$ poly $(LOPI_1)(LOPI_1)$ (constant in equation throughout convolution, no growth) therefore

$Constant1 + Const(b) + (Const(b-1) + Const2(b^2-b)$

In all the polynomials constants excluded as they do not change O time calculation, the variables show $O(x^3)$ for $LOPI^2$ subsets when $LOPI_1=LOPI_2$

Adapt a=b in subsets of Superset LOPI 17, as all twin LOPI are used to generate NP elements in each subset, to generate possible secondary twin prime products in subsets of LOPI 17, **set a=b and query for multiplicity 1 in subsets in range as discussed.**

f(a,b), a=b, a=b+1, for [twin sub 5,7], [twin sub 11,13] but for [twin sub 17,1] poly is

$C + C_1(b+1) + C_2(b) + C_3(b^2 + b)$  therefore becomes $O(x^3)$  [15 pg.107-110]

**3**-Query 3 subset matrices for the LOPI above of the Integer in grid area LOPI r min – LOPI max cycle= 5 min, Max= number of cycles of 18≈ (2267/5 ) (+ 17)

Query x=2267 LOPI 17

grid area m=5, n=25   subsets LOPI 17 {B,C,D}

Stop Output : Element not found

0 = Prime

New Query : 6313  LOPI 13

Query x = 6313

grid area m=0, n=70 subsets LOPI 13 {B,C,D,E}

Output : Found at sub B (4, 6)

Element not found sub C, D, E,

Stop Output

1 find = Secondary

Output : Found : sub B (4, 6)   (m-1 = 3, n-1 = 5) = a, b

(Put into linear Dirichlet binomials in convolution equation for subB  LOPI 13

Convolution equation LOPI 13 sub B : equation : $(5 + 18(a))(17 + 18(b))$

 Function $(5 + 18(3))(17 + 18(5))$

Output : 59    107    (factors of Query 6313)

Check 59 * 107 = 6313

Stop Output

Discussion:

By using the least LOPI residue m=LOPI 5 as a divisor for Query Integer and the result divided by 18 we get the maximum number of cycles of 18 in which to search for the integer,

Integer/5 = cycles of 18,=a in $(17+ 18(a))$ max range.

Ex :889=  search  grid LOPI 7 m=5, n=10   we find factors between m=5 and n=10   between 5-179. Factors are 7 and 127.

N=10 is based Integer/5, 889/5=about 179/18 which about 10 cycles of 18, therefore 18(10)+LOPImax -17 = 197, which is the area of grid to search possible factors m-5  to n-10 .

Based on this information we can target the area of our query thereby reducing query string length in the data function table to optimize time efficiency.  We do not have to search the full length of the matrix, just the areas specific to the task.

<u>Adaptation to generate possible set of secondary twin prime products to LOPI class 17 subsets:</u>

Adapt a=b in subsets of Superset LOPI 17, as all twin LOPI are used to generate NP elements in each subset, to generate possible secondary twin prime products in subsets of LOPI 17**, set a=b and query for multiplicity 1 in subsets in range as**

**discussed.**

f(a,b),a=b, a=b+1, for [twin sub 5,7] , [twin sub 11,13]

**$LOPI_1(LOPI_2)$** + **$(LOPI_1(18))$**(a) + **$(LOPI_2 18)$**(a) + **$18^2$**($a^2$)  [15 pg.107-110]

$LOPI_1$, $LOPI_2$ are the twin LOPI pairs, [5,7] [11,13] that create the constants:

above, and $18^2$ never grows either, as it is 324, the last term in the general functions above where a≥0, b≥0 therefore

$18^2(ab)=324(ab)$, C(ab).

but for [twin sub 17,1] poly is

$17 + 306(b+1) + 18(b) + 324(b^2 + b)$   $O(x^3)$  [Ref. 1 pg.107-110]

The integers in the equations are constants as they are made on the constant LOPI identity, only a and b change in the LOPI Dirichlet Linear Convolutions of LOPI Mod 18.

**Goldbach/Dirichlet LOPI Congruence Seesaw Algorithm using ELMA-18**

Question: solutions for 2N>=6, with Prime + Prime = 2N

-Sum LEPI

-Calculate the linear progressions for a=0 or 1 for LOPI 1, to b=c, c-1, b decided by c, cycles of 18 of the LEPI. Even Integer – LEPI/18 = c = b or c-1 = b. We are subtracting $b_{max}$ by 1 each time and adding 1 each time to $a_0$. O(x) [2]

LEPI + 18(c) ≡ $LOPI_1$ + 18(a) (Mod $LOPI_2$ + 18(b)) with bmax i=c, $a_i$=0, depending on LEPI and specific LOPIs needed from template information in article.

LEPI + 18(c) ≡ $LOPI_1$ + 18(a) (Mod $LOPI_2$ + 18(b)), progressions $a_i$=0, $bmax_i$=c,c-1

Multiply the results from this area of each linear progression in the residue to the corresponding result in the modulo linear progression. a=0 to [$a_{f\,max}$=$b_{i\,max}$] and $bmax_i$ to [b=0=$a_i$] as shown in the model of the Goldbach/Dirichlet LOPI Congruence Seesaw for integer 88, from page 58 in the article,

**88** (16+18(c))≡(5 + 18(a) (Mod 11 + 18(b)  $a_{min}$=0, $b_{max}$=c c=4    LEPI 16 / 5+(11) e (11*5) LOPI I

| | | |
|---|---|---|
| 5    $a_i$=0 | 83    $b_{i=c}$4 | 5*83=415 |
| 23 | 65 | 23*65=1495 |
| 41 | 47 | 41*47=1927 |
| 59 | 29 | 59*29=1711 |
| 77 $a_f$ max, | 11    $b_f$=0 | 77*11=847 |

At this point we can choose to query each element of the addition operation for multiplicity in their respective obvious subsets or perhaps more efficient method is to multiply the two elements that add to give the even integer and query their product multiplicity in the subset that we already know, LOPI 1 Superset / Subset 11(5).

We can query in the range needed for the maximum value of cycles needed as shown above in the algorithm query section for Primality testing and finding factors.

Query non-prime subsets of this LOPI class for multiplicity

grid area range as min a=0, to max product treated as above, Max/5 cycles(18)+17

1= secondary = prime + prime Goldbach solution for Integer the question 2N≥6,

Data Even Integer:  LEPI  Look at table for predetermined template pg. 58

Calculate See saw congruence LOPI  relations

Query each element in each pair of congruences that are generated

Or

Query known LOPI superset's subset for multiplicity of product of potential Goldbach Prime pairs  generated.

For example in the above case we know 88=LEPI 16 made by LOPI 11 + LOPI 5 in the template used, although we haven't finished producing all possible P+P=88 until we make the seesaw linear progressions for each template.  LOPI 16 is also made by LOPI 17 + LOPI 17, but for this example we provided only the seesaw solutions for LOPI 5 + LOPI 11.

The product of these two LOPIs will be 55 which in in the LOPI 1 reduced residue class mod 18 in our work.  Also the template used shows us the subset to query, 11*5 seed subset in LOPI 1.

Upon our query in the appropriate grid area we will get a multiplicity of 1 for all secondary non-primes in the LOPI 1 subset of (11*5), or a P * P, which are our Goldbach prime elements that add to equal 2N≥6.

In this simple case we know we have already 3 solutions in such a small integer and haven't yet included the solutions to the LOPI template of 17+17, which is one 71*17.

Multiplicity 1=  output found once

Products with Mult=1        a) 1927=41*47    415=5*83,         1711 = 59*29

Secondary primes that are produced by the multiplication of two primes to yield the even LEPI integer 88 when added.